\documentclass[11pt]{article}

\usepackage[utf8]{inputenc}
\usepackage{lmodern}
\usepackage[T1]{fontenc}

\usepackage{amsmath, amssymb, amsfonts, amsthm}
\usepackage{mathtools}
\usepackage{bm}
\usepackage{dsfont}
\usepackage{latexsym}
\usepackage{cases}

\usepackage{caption}
\usepackage{subcaption}

\usepackage{graphicx}
\usepackage{float}
\usepackage{tikz, pgfplots}
\usetikzlibrary{3d, hobby}
\pgfplotsset{compat=1.18}

\usepackage{tabularx}
\usepackage{booktabs}
\usepackage{multirow}
\usepackage{array}
\usepackage{cellspace}
\renewcommand{\arraystretch}{1.5}

\usepackage{nicefrac}
\usepackage{microtype}
\usepackage{color}
\usepackage{setspace}
\usepackage{verbatim}

\usepackage{enumerate}
\usepackage{xpatch}
\usepackage{enumitem}

\usepackage[hyphens]{url}
\usepackage[colorlinks=true,allcolors=blue,breaklinks=true]{hyperref}
\usepackage[authoryear]{natbib}

\usepackage[a4paper, left=1in, right=1in, top=1in, bottom=1in]{geometry}

\newtheorem{assumption}{Regularity Conditions}[section]

\newtheorem{theorem}{Theorem}[section]
\newtheorem{remark}[theorem]{Remark}
\newtheorem{lemma}[theorem]{Lemma}

\newtheorem{corollary}[theorem]{Corollary}
\numberwithin{equation}{section}
\allowdisplaybreaks

\title{\textbf{Uniform Asymptotic Theory for Local Likelihood Estimation of Covariate-Dependent Copula Parameters}}

\author{ %my God the father and Jesus Christ and
	\small Mathias N. Muia$^{\dagger}$\\[0.5em]
	\scriptsize $^{\dagger}$Department of Mathematics and Statistics, University of South Alabama, Mobile, AL 36688, USA\\
	[0.5em]
	\scriptsize $^{\dagger}$\texttt{mnmuia@southalabama.edu} 
}
\date{}
\begin{document}
	\maketitle
	
\begin{abstract}
	Conditional copula models allow dependence structures to vary with observed covariates while preserving a separation between marginal behavior and association. We study the uniform asymptotic behavior of kernel-weighted local likelihood estimators for smoothly varying copula parameters in multivariate conditional copula models. Using a local polynomial approximation of a suitably transformed calibration function, we establish uniform convergence rates over compact covariate sets for the local log-likelihood, its score, and its Hessian. These results yield uniform consistency of the local maximum likelihood estimator and of the induced copula parameter function. The analysis is based on empirical process techniques for kernel-indexed classes with shrinking neighborhoods and polynomial entropy bounds, providing theoretical support for global consistency and stable local optimization in covariate-dependent copula models.
\end{abstract}


\section{Introduction}

Copulas offer a flexible framework for modeling dependence structures independently of marginal distributions. By Sklar’s theorem~\cite{sklar1959fonctions}, any multivariate distribution with continuous marginals admits a unique copula representation, which allows dependence to be studied separately from marginal behavior. In many contemporary applications, however, dependence is not constant and often varies with observable covariates. Such covariate-driven dependence is common in fields such as finance, environmental science, and biomedical research, where associations among multiple outcomes evolve across time, space, or demographic characteristics.

Conditional copula models are designed to capture this heterogeneity by extending Sklar’s representation to conditional distributions~\cite{patton2006modelling}. Let $(\mathbf{X},\mathbf{Y})=(X_1,\ldots,X_n,Y_1,\ldots,Y_s)$ be a random vector, and define the conditional marginal distributions by
$F_{k\mid \mathbf{Y}}(x_k\mid \mathbf{Y}=\mathbf{y})=\mathbb{P}(X_k\le x_k\mid \mathbf{Y}=\mathbf{y})$, $k=1,\ldots,n$.
The conditional joint distribution of $\mathbf{X}$ given $\mathbf{Y}=\mathbf{y}$ is
\(
H_{\mathbf{X}\mid \mathbf{Y}}(x_1,\ldots,x_n\mid \mathbf{y})
=
\mathbb{P}(X_1\le x_1,\ldots,X_n\le x_n\mid \mathbf{Y}=\mathbf{y}).
\)
When the conditional marginals are continuous, there exists a unique $n$-dimensional conditional copula $C(\cdot\mid \mathbf{y})$ satisfying
\(
H_{\mathbf{X}\mid \mathbf{Y}}(x_1,\ldots,x_n\mid \mathbf{y})
=
C\!\left(
F_{1\mid \mathbf{Y}}(x_1\mid \mathbf{y}),\ldots,
F_{n\mid \mathbf{Y}}(x_n\mid \mathbf{y})
\mid \mathbf{y}
\right).
\)
Together, the conditional copula and the conditional marginals fully characterize the conditional joint distribution.

This formulation implies that the conditional pseudo-observations
$U_{k\mid \mathbf{Y}}=F_{k\mid \mathbf{Y}}(X_k\mid \mathbf{Y})$, $k=1,\ldots,n$, are uniformly distributed on $(0,1)$ and satisfy
\(
(U_{1\mid \mathbf{Y}},\ldots,U_{n\mid \mathbf{Y}})\mid \mathbf{Y}=\mathbf{y}
\sim C(\cdot\mid \mathbf{y}).
\)
Consequently, the conditional copula captures the entire dependence structure of $\mathbf{X}$ at a given covariate value $\mathbf{Y}=\mathbf{y}$. In parametric conditional copula models, covariate effects are typically encoded through a smooth parameter function $\theta(\mathbf{y})$, leading to copulas of the form $C(\cdot\mid \theta(\mathbf{y}))$. Common examples include Archimedean families with parametric generators~\cite{mcneil2009multivariate,genest1995semiparametric}.

Inference in this setting requires estimation of the unknown calibration function $\theta(\cdot)$ while respecting the constraints imposed by the copula family. Local likelihood methods provide a natural nonparametric strategy for this task. The approach consists of working with a transformed calibration function $\nu(\mathbf{y})=\psi(\theta(\mathbf{y}))$, approximating it locally by multivariate polynomials, and maximizing a kernel-weighted conditional log-likelihood. This construction yields estimators that combine flexibility with favorable asymptotic behavior.

A substantial body of work addresses estimation in copula models under various assumptions. Kernel-based estimation of copula densities and related quantities, together with uniform consistency and asymptotic normality, has been studied in early contributions~\cite{gijbels1990estimating}. Subsequent work has focused on mitigating boundary bias and developing transformation-based corrections, leading to improved weak convergence properties~\cite{omelka2009improved}. Semiparametric estimation of copula parameters without specifying marginal distributions has also been developed using minimum distance and rank-based methods~\cite{tsukahara2005semiparametric}.

Research on covariate-dependent dependence structures has expanded considerably in recent years. Semiparametric conditional copula models with unknown calibration functions have been proposed and estimated via local polynomial pseudo-likelihood methods~\cite{abegaz2012semiparametric}. Local likelihood approaches for parametric conditional copulas with covariate-varying parameters have led to pointwise asymptotic bias and variance expressions, along with data-driven model selection and inference procedures~\cite{acar2011dependence}. Fully nonparametric estimators of conditional copulas and associated measures of conditional dependence have also been investigated~\cite{veraverbeke2011estimation}.

Related developments include dynamic copula models with time-varying dependence parameters estimated through local likelihood, which have found applications in financial time series analysis~\cite{hafner2010efficient}. Comparisons of parametric, semiparametric, and inference-function-for-margins approaches highlight the robustness of semiparametric methods under marginal misspecification~\cite{kim2007comparison}. More recent work has examined nonparametric estimation of copula densities, with particular attention to boundary effects and unbounded behavior near the corners, and has proposed kernel-based procedures supported by theoretical guarantees and empirical evidence~\cite{muia2025kernel,muia2025kernel1}.

Despite these advances, most existing asymptotic results for covariate-dependent copula models are pointwise, describing behavior at fixed covariate values. Such results do not address several issues of theoretical and practical importance. Uniform control over the covariate space is required to establish global consistency, to ensure stability of numerical optimization procedures, and to justify data-driven bandwidth selection. Uniform convergence of the score and Hessian processes also plays a key role in analyzing local maximizers and deriving higher-level asymptotic properties.

The present paper addresses these gaps by developing a uniform asymptotic theory for local likelihood estimators in covariate-dependent copula models. While earlier work provides pointwise bias and variance expansions~\cite{acar2011dependence,muia2025multivariate}, we establish uniform convergence of the kernel-weighted local log-likelihood, its gradient, and its Hessian over compact subsets of the covariate space. Achieving this requires careful control of empirical processes indexed jointly by covariate locations and local polynomial coefficients, as well as handling the shrinking supports induced by kernel localization.

Our main objective is to prove uniform consistency for local likelihood estimators of covariate-dependent copula parameters. We derive explicit uniform convergence rates for the local log-likelihood, the score, and the Hessian, which in turn imply uniform consistency of the local maximum likelihood estimator of the calibration function $\nu(\mathbf{y})$ and the induced copula parameter function $\theta(\mathbf{y})$. Although these results provide theoretical support for data-driven bandwidth selection, a detailed asymptotic analysis of the selected bandwidth itself is beyond the scope of this paper.

The analysis relies on empirical process techniques for classes of functions indexed by both covariate locations and local polynomial coefficients. The main difficulty stems from multivariate kernel localization, which produces shrinking neighborhoods and diverging envelopes as the bandwidth decreases. These challenges are addressed using geometric properties of the covariate space, polynomial entropy bounds, and maximal inequalities for empirical processes. As is typical in nonparametric smoothing, the resulting convergence rates depend on the covariate dimension $s$, reflecting the familiar curse of dimensionality.

The contributions of this paper are threefold. First, we establish uniform stochastic equicontinuity for kernel-weighted local likelihood processes under mild regularity conditions. Second, we obtain explicit uniform convergence rates that capture the effects of bandwidth choice, sample size, and covariate dimension. Third, we provide a theoretical foundation for uniform inference and reliable global optimization in covariate-dependent copula models.

The remainder of the paper is organized as follows. Section~\ref{sec:setup} introduces the conditional copula model and the local likelihood estimation framework. Section~\ref{sec:tuning} discusses implementation issues, including bandwidth selection and copula family selection using cross-validation criteria. Section~\ref{sec:main-results} presents the main theoretical results, including uniform convergence rates for the local log-likelihood, its score, and its Hessian, and the resulting consistency statements. Section~\ref{sec:simulation} reports simulation results illustrating finite-sample performance. The Appendix contains additional background material and technical arguments.

\section{Model Setup and Local Likelihood Estimation}
\label{sec:setup}

Let $(\mathbf{X},\mathbf{Y})=(X_1,\ldots,X_n,Y_1,\ldots,Y_s)$ be a random vector, where $\mathbf{Y}\in\mathbb{R}^s$.
Consider a parametric copula family $\{C(\cdot\mid\theta):\theta\in\Theta\}$ with associated copula density $c(\cdot\mid\theta)$.
We assume that, conditional on $\mathbf{Y}=\mathbf{y}$, the random vector $\mathbf{X}$ has copula $C(\cdot\mid\theta(\mathbf{y}))$, where the copula parameter $\theta(\mathbf{y})$ is an unknown smooth function taking values in a parameter space $\Theta\subset\mathbb{R}^d$.
Our objective is to estimate the calibration function $\theta(\cdot)$ by locally approximating it with a multivariate polynomial of degree $p$.

Since the copula parameter function $\theta(\mathbf{y})$ must take values in the admissible parameter space $\Theta$, we introduce a strictly monotone link transformation
\(
\nu(\mathbf{y})=\psi(\theta(\mathbf{y})), \qquad
\theta(\mathbf{y})=\psi^{-1}(\nu(\mathbf{y})),
\)
where the inverse link $\psi^{-1}:\mathbb{R}\to\Theta$ maps the real line onto the parameter space, ensuring that $\theta(\mathbf{y})\in\Theta$ for all $\mathbf{y}$.
The transformed function $\nu(\mathbf{y})$ is referred to as the \emph{calibration function}.
This transformation is essential because local polynomial modeling is naturally formulated on $\mathbb{R}$, whereas copula parameters are typically subject to constraints.
For example, the Gaussian copula parameter satisfies $\theta\in(-1,1)$, while for the Clayton copula one commonly assumes $\theta\in(0,\infty)$.
The link function enforces these constraints throughout the estimation procedure.

We next describe the local polynomial approximation of the calibration function.
Let $\{(X_{1i},\ldots,X_{ni},Y_{1i},\ldots,Y_{si})\}_{i=1}^N$ be a random sample, and suppose that the calibration function $\nu(\cdot)$ is $(p+1)$-times continuously differentiable at an interior point $\mathbf{y}\in\mathbb{R}^s$.
A multivariate Taylor expansion of $\nu(\mathbf{Y}_i)$ about $\mathbf{y}$ yields
\begin{align}\label{multi-taylor}
	\nu(\mathbf{Y}_i)
	&= \nu(\mathbf{y})
	+ \nabla \nu(\mathbf{y})^\top (\mathbf{Y}_i - \mathbf{y})
	+ \frac{1}{2} (\mathbf{Y}_i - \mathbf{y})^\top \nabla^2 \nu(\mathbf{y}) (\mathbf{Y}_i - \mathbf{y})
	+ \cdots
	+ R_{p+1}(\mathbf{Y}_i - \mathbf{y}),
\end{align}
where $\nabla \nu(\mathbf{y})$ and $\nabla^2 \nu(\mathbf{y})$ denote the gradient and Hessian of $\nu$ evaluated at $\mathbf{y}$, respectively, and the remainder term satisfies $R_{p+1}(\mathbf{Y}_i - \mathbf{y})=O(\|\mathbf{Y}_i-\mathbf{y}\|^{p+1})$.
For notational convenience, we approximate $\nu(\mathbf{Y}_i)$ by
\(
\nu(\mathbf{Y}_i)\approx \boldsymbol{\phi}_p(\mathbf{Y}_i-\mathbf{y})^\top \boldsymbol{\gamma},
\)
where $\boldsymbol{\phi}_p(\cdot)$ collects all monomials of total degree at most $p$, and $\boldsymbol{\gamma}$ denotes the corresponding vector of local polynomial coefficients.

Based on this approximation, we construct a kernel-weighted local log-likelihood.
Let $\mathbf{U}_i = (U_{1i}, \ldots, U_{ni})^\top$ denote pseudo-observations from the conditional marginals, with
\(
U_{ki} = F_{k \mid \mathbf{Y}}(X_{ki} \mid \mathbf{Y}_i).
\)
Assuming conditional independence across observations, the contribution of each data point $(\mathbf{U}_i \mid \mathbf{Y}_i)$ in a neighborhood of $\mathbf{y}$ to the log-likelihood is given by
\(
\log c(\mathbf{U}_i \mid \psi^{-1}(\nu(\mathbf{Y}_i))),
\)
where $c(\cdot\mid\cdot)$ denotes the conditional copula density function.
The resulting kernel-weighted local log-likelihood is
\begin{align}\label{multi-local-lik}
	L_N(\boldsymbol{\gamma},\mathbf{y},p,h)
	=
	\sum_{i=1}^{N}
	\log c\!\left(
	\mathbf{U}_i \,\middle|\,
	\psi^{-1}\!\bigl(\boldsymbol{\phi}_p(\mathbf{Y}_i-\mathbf{y})^\top\boldsymbol{\gamma}\bigr)
	\right)
	K_h(\mathbf{Y}_i-\mathbf{y}),
\end{align}
where $K_h(\mathbf{v})=h^{-s}K(\mathbf{v}/h)$ is a rescaled multivariate kernel.

The kernel weight $K_h(\mathbf{Y}_i-\mathbf{y})$ localizes estimation to observations with covariates near $\mathbf{y}$.
More generally, one may consider a multivariate kernel with bandwidth matrix $H$ of the form $K_H(\mathbf{v})=|H|^{-1}K(H^{-1}\mathbf{v})$; see, e.g.,~\cite{wand1995kernel}.
For simplicity, we adopt the isotropic specification $H=hI_s$, which applies uniform smoothing across all covariate dimensions while preserving standard normalization and asymptotic properties.
Typical choices for $K$ include the product Epanechnikov kernel and the multivariate Gaussian kernel.

To characterize the maximizer of~\eqref{multi-local-lik}, we define the score vector and Hessian matrix with respect to the local coefficient vector $\boldsymbol{\gamma}$.
Let $\ell(\theta,\mathbf{u})=\log c(\mathbf{u}\mid\theta)$ denote the log-copula density, and write $(\cdot)=(\boldsymbol{\phi}_p(\mathbf{Y}_i-\mathbf{y})^\top\boldsymbol{\gamma})$.
Applying the chain rule yields the score function
\begin{align*}
	\nabla_{\boldsymbol{\gamma}}L_N(\boldsymbol{\gamma},\mathbf{y})
	=
	\sum_{i=1}^N
	\ell'\!\left(\psi^{-1}(\cdot),\mathbf{U}_i\right)
	(\psi^{-1})'(\cdot)\,
	\boldsymbol{\phi}_p(\mathbf{Y}_i-\mathbf{y})\,
	K_h(\mathbf{Y}_i-\mathbf{y}),
\end{align*}
and the Hessian matrix
\begin{align*}
	\nabla^2_{\boldsymbol{\gamma}}L_N(\boldsymbol{\gamma},\mathbf{y})
	=
	\sum_{i=1}^N
	\Big[
	\ell''\!\left(\psi^{-1}(\cdot),\mathbf{U}_i\right)\{(\psi^{-1})'(\cdot)\}^2
	+
	\ell'\!\left(\psi^{-1}(\cdot),\mathbf{U}_i\right)(\psi^{-1})''(\cdot)
	\Big]
	\\[-2mm]
	\hspace{3cm}\times
	\boldsymbol{\phi}_p(\mathbf{Y}_i-\mathbf{y})
	\boldsymbol{\phi}_p(\mathbf{Y}_i-\mathbf{y})^\top
	K_h(\mathbf{Y}_i-\mathbf{y}).
\end{align*}

The local maximum likelihood estimator is defined as
\begin{align}\label{estimator_main}
	\hat{\boldsymbol{\gamma}}(\mathbf{y})
	=
	\arg\max_{\boldsymbol{\gamma}}
	L_N(\boldsymbol{\gamma},\mathbf{y},p,h),
\end{align}
and is obtained by solving the estimating equation $\nabla_{\boldsymbol{\gamma}}L_N(\boldsymbol{\gamma},\mathbf{y})=0$.
In practice, this maximization can be carried out using a Newton--Raphson algorithm, with updates of the form
\begin{align}\label{newton-raphson}
	\boldsymbol{\gamma}^{(m+1)}
	=
	\boldsymbol{\gamma}^{(m)}
	-
	\left\{
	\nabla^2_{\boldsymbol{\gamma}}L_N(\boldsymbol{\gamma}^{(m)},\mathbf{y})
	\right\}^{-1}
	\nabla_{\boldsymbol{\gamma}}L_N(\boldsymbol{\gamma}^{(m)},\mathbf{y}),
\end{align}
iterated until convergence.

The resulting estimator of the calibration function at $\mathbf{y}$ is given by
$\hat{\nu}(\mathbf{y})=\mathbf{e}_0^\top\hat{\boldsymbol{\gamma}}(\mathbf{y})$, where $\mathbf{e}_\alpha$ denotes the unit vector selecting the intercept ($\alpha^{\text{th}}$-order) coefficient in the local polynomial expansion.
The corresponding estimator of the copula parameter is then
$\hat{\theta}(\mathbf{y})=\psi^{-1}(\hat{\nu}(\mathbf{y}))$.
Substituting $\hat{\theta}(\mathbf{y})$ into $C(\cdot\mid\theta(\mathbf{y}))$ yields an estimator of the copula function at the covariate value $\mathbf{y}$.

It is important to note that the above maximization problem is localized at a fixed evaluation point $\mathbf{y}$ and therefore yields an estimate of the copula parameter only at that point.
Consequently, to recover the full calibration function $\theta(\cdot)$, the local maximization problem must be solved over a sufficiently fine grid of points $\mathbf{y}$ spanning the domain of the covariate $\mathbf{Y}$.
For completeness, additional asymptotic properties of the local likelihood estimator, including pointwise bias and variance expansions as well as asymptotic normality, are discussed in the Appendix.

The primary focus of this paper is to establish the uniform asymptotic behavior of these estimators as the evaluation point $\mathbf{y}$ varies over the compact set $\mathcal{Y}_0$.

\section{Model Tuning}
\label{sec:tuning}

Tuning a covariate-dependent copula model involves two related decisions. The first concerns the choice of the smoothing bandwidth used in local likelihood estimation. The second concerns selection of an appropriate parametric copula family. We adopt a cross-validation strategy that is standard in local likelihood inference and follows the general ideas developed in~\cite{acar2011dependence}, adapted here to multivariate covariates.

\subsection{Bandwidth Selection}

Because estimation is based on a kernel-weighted local copula likelihood, bandwidth selection is naturally handled by leave-one-out cross-validation. Let $\hat{\theta}_h(\cdot)$ denote the estimator of the copula parameter function obtained with bandwidth $h$. For each observation $(\mathbf{U}_i,\mathbf{Y}_i)$, a leave-one-out estimate $\hat{\theta}_{h,-i}(\mathbf{Y}_i)$ is computed using all observations except the $i$th. The bandwidth is then chosen by maximizing the cross-validated local likelihood
\begin{align}\label{lcv_h}
	\mathrm{CVL}(h)
	=
	\sum_{i=1}^N
	\log c\!\left(\mathbf{U}_i \mid \hat{\theta}_{h,-i}(\mathbf{Y}_i)\right).
\end{align}
The selected bandwidth $h_{cv}$ is the maximizer of~\eqref{lcv_h}. In settings with highly uneven covariate designs or strong local variation, adaptive choices such as variable or nearest-neighbor bandwidths may also be employed.

\subsection{Copula Family Selection}

Choosing the bandwidth does not determine the form of dependence, which is governed by the copula family. Since likelihood values are not directly comparable across different families, selection is based on predictive performance rather than likelihood maximization.

Let $\mathcal{C}=\{C_q:q=1,\ldots,Q\}$ denote a collection of candidate copula families. For each family $C_q$, the bandwidth selection procedure yields an optimal bandwidth $h_q^\ast$, which is then used to compute leave-one-out estimates $\hat{\theta}^{(-i)}_{h_q^\ast}(\mathbf{Y}_i)$. These estimates define candidate conditional copula models of the form $C_q(\cdot\mid\hat{\theta}^{(-i)}_{h_q^\ast}(\mathbf{Y}_i))$.

Predictive accuracy is evaluated through conditional prediction of each pseudo-observation given the remaining components and the covariate. For a given family $C_q$, the conditional expectation of $U_{ki}$ given $\mathbf{U}_{-k,i}$ and $\mathbf{Y}_i$ is
\[
\widehat{\mathbb{E}}^{(-i)}_q
\!\left(U_{ki}\mid \mathbf{U}_{-k,i},\mathbf{Y}_i\right)
=
\int_0^1
u_k\,
c_{q,k}\!\left(
u_k \mid \mathbf{U}_{-k,i},
\hat{\theta}^{(-i)}_{h_q^\ast}(\mathbf{Y}_i)
\right)\,du_k,
\]
where $c_{q,k}(\cdot\mid\cdot)$ denotes the conditional copula density implied by family $C_q$.

The corresponding cross-validated prediction error is defined as
\[
\mathrm{CVPE}(C_q)
=
\sum_{i=1}^N
\sum_{k=1}^n
\left\{
U_{ki}
-
\widehat{\mathbb{E}}^{(-i)}_q
\!\left(U_{ki}\mid \mathbf{U}_{-k,i},\mathbf{Y}_i\right)
\right\}^2.
\]
The copula family that minimizes $\mathrm{CVPE}(C_q)$ over $q=1,\ldots,Q$ is selected.

This criterion can be interpreted as an empirical approximation of a conditional mean squared prediction error. It is minimized when the working copula family matches the underlying dependence structure, providing a likelihood-free and interpretable basis for copula family selection.

\section{Main Results}\label{sec:main-results}

We begin by introducing additional notation and specifying the probabilistic setting in which the main results are derived. Let $\mathcal{Y}_0 \subset \mathbb{R}^s$ be a nonempty compact set contained in the interior of the support of the covariate vector $\mathbf{Y}$. We assume that $\mathbf{Y}$ admits a density $f_{\mathbf{Y}}$ with respect to Lebesgue measure. All uniform statements in this section are understood to hold over $\mathcal{Y}_0$.

Throughout, we work with the scaled local log-likelihood
\[
\tilde L_N(\boldsymbol{\gamma};\mathbf{y})
=
(Nh^s)^{-1} L_N(\boldsymbol{\gamma};\mathbf{y}),
\]
and its population counterpart
\[
Q(\boldsymbol{\gamma};\mathbf{y})
=
\mathbb{E}\!\left[\tilde L_N(\boldsymbol{\gamma};\mathbf{y})\right].
\]
The results below describe the uniform behavior of the local log-likelihood, its score, and its Hessian as stochastic processes indexed jointly by the evaluation point $\mathbf{y}$ and the local polynomial coefficient vector $\boldsymbol{\gamma}$.

\begin{assumption}\label{assump:regularity}
	The following conditions are imposed.
	\begin{enumerate}
		\item[(R1)]
		The covariate vector $\mathbf{Y} \in \mathbb{R}^s$ has a density $f_{\mathbf{Y}}$ with respect to Lebesgue measure. There exist constants $0 < c_f \le C_f < \infty$ such that
		$c_f \le f_{\mathbf{Y}}(\mathbf{y}) \le C_f$ for all $\mathbf{y} \in \mathcal{Y}_0$, and $f_{\mathbf{Y}}$ is continuous on an open neighborhood of $\mathcal{Y}_0$.
		
		\item[(R2)]
		The kernel $K : \mathbb{R}^s \to \mathbb{R}$ is bounded, symmetric, Lipschitz continuous, and compactly supported on $[-1,1]^s$. It satisfies $\int_{\mathbb{R}^s} K(v)\,dv = 1$ and has finite $(p+2)$-moment
		$\int_{\mathbb{R}^s} \|v\|^{p+2} |K(v)|\,dv < \infty$.
		
		\item[(R3)]
		The bandwidth sequence $h=h_N$ satisfies $h_N \to 0$ as $N \to \infty$, $N h_N^s / \log(1/h_N) \to \infty$, and $N h_N^{s+2(p+1)} \to 0$.
		
		\item[(R4)]
		The calibration function $\nu : \mathbb{R}^s \to \mathbb{R}$ is $(p+2)$-times continuously differentiable on an open neighborhood of $\mathcal{Y}_0$, with all partial derivatives up to total order $p+2$ uniformly bounded.
		
		\item[(R5)]
		For each $\theta$ in a compact set $\Theta_0$ containing $\{\theta(\mathbf{y}) : \mathbf{y} \in \mathcal{Y}_0\}$, the copula log-density $\ell(\theta,u)=\log c(u\mid\theta)$ is twice continuously differentiable in $\theta$ for all $u \in (0,1)^n$. There exist envelope functions $M_1$ and $M_2$ such that
		\[
		\sup_{\theta\in\Theta_0} |\ell'(\theta,u)| \le M_1(u),
		\qquad
		\sup_{\theta\in\Theta_0} |\ell''(\theta,u)| \le M_2(u),
		\]
		with $E[M_1(U)^2]<\infty$ and $E[M_2(U)^2]<\infty$. The conditional Fisher curvature
		\[
		\sigma^2(\mathbf{y})
		=
		-
		E\!\left[\ell''(\theta(\mathbf{y}),U)\mid \mathbf{Y}=\mathbf{y}\right]
		\]
		is continuous in $\mathbf{y}$ and bounded away from zero and infinity on $\mathcal{Y}_0$.
		
		\item[(R6)]
		The inverse link $\psi^{-1}$ is strictly monotone and twice continuously differentiable, with bounded first and second derivatives on $\{\nu(\mathbf{y}): \mathbf{y}\in\mathcal{Y}_0\}$.
		
		\item[(R7)]
		For each $\mathbf{y}\in\mathcal{Y}_0$, the population criterion $Q(\boldsymbol{\gamma};\mathbf{y})$ admits a unique maximizer $\boldsymbol{\gamma}^\ast(\mathbf{y})$ in a compact set $\Gamma \subset \mathbb{R}^{d_p}$. The map $\mathbf{y}\mapsto\boldsymbol{\gamma}^\ast(\mathbf{y})$ is continuous, and the population Hessian
		\[
		H(\boldsymbol{\gamma};\mathbf{y})
		=
		E\!\left[
		\nabla^2_{\boldsymbol{\gamma}}
		\ell\!\left(
		\psi^{-1}(\phi_p(\mathbf{Y}-\mathbf{y})^\top \boldsymbol{\gamma}),U
		\right)
		K\!\left(\tfrac{\mathbf{Y}-\mathbf{y}}{h}\right)
		\right]
		\]
		is uniformly negative definite on $\Gamma\times\mathcal{Y}_0$.
		
		\item[(R8)]
		There exist constants $C_\ell^{(0)}, C_\ell^{(1)}, C_\ell^{(2)} < \infty$ such that
		\[
		\sup_{\theta\in\Theta_0}\sup_{u\in(0,1)^n}
		\big|\ell^{(\kappa)}(\theta,u)\big|
		\le C_\ell^{(\kappa)}, \qquad \kappa=0,1,2.
		\]
	\end{enumerate}
\end{assumption}

\begin{lemma}\label{lem:kernel-class-entropy}
	Assume \textnormal{(R1)}, \textnormal{(R2)}, \textnormal{(R4)}, \textnormal{(R5)}, 
	\textnormal{(R6)}, and \textnormal{(R8)}, and let $\mathcal{Y}_0\subset\mathbb{R}^s$ and
	$\Gamma\subset\mathbb{R}^{d_p}$ be compact. For each
	$(\mathbf{y},\boldsymbol{\gamma})\in\mathcal{Y}_0\times\Gamma$, define
	\[
	f_{\mathbf{y},\boldsymbol{\gamma}}(U,\mathbf{Y})
	=
	h^{-s}\,
	\ell\!\Big(
	\psi^{-1}(\phi_p(\mathbf{Y}-\mathbf{y})^\top \boldsymbol{\gamma}),\,U
	\Big)\,
	K\!\left(\frac{\mathbf{Y}-\mathbf{y}}{h}\right),
	\]
	and let
	\(
	\mathcal{F}_h
	=
	\{ f_{\mathbf{y},\boldsymbol{\gamma}} : \mathbf{y}\in\mathcal{Y}_0,\ 
	\boldsymbol{\gamma}\in\Gamma \}.
	\)
	Then, as $h\to0$, the following statements hold.
	\begin{enumerate}[label=(\roman*)]
		\item There exists an envelope $F_h$ for $\mathcal{F}_h$ and constants
		$0<c_1\le c_2<\infty$, independent of $h$, such that
		\[
		c_1\,h^{-s}
		\;\le\;
		\|F_h\|_{L_2(P)}
		\;\le\;
		c_2\,h^{-s}.
		\]
		
		\item There exists a constant $C>0$, independent of $h$, such that for all
		$0<\varepsilon\le1$,
		\[
		\log N\!\left(
		\varepsilon \|F_h\|_{L_2(P)},\
		\mathcal{F}_h,\
		L_2(P)
		\right)
		\le
		C\!\left\{1+\log(1/h)+\log(1/\varepsilon)\right\}.
		\]
		
		\item Let $\mathbb{P}_N$ denote the empirical measure based on
		$\{(U_i,\mathbf{Y}_i)\}_{i=1}^N$. Then
		\[
		\sup_{f\in\mathcal{F}_h}
		\big|(\mathbb{P}_N-P)f\big|
		=
		O_p\!\left(
		\sqrt{\frac{\log(1/h)}{N h^s}}
		\right).
		\]
	\end{enumerate}
\end{lemma}

\begin{remark}\label{rem:logN-rate}
	The logarithmic factor in Lemma~\ref{lem:kernel-class-entropy} arises from the entropy bound. Under the mild additional condition
	\[
	\log(1/h)=O(\log N),
	\]
	the stochastic rate may equivalently be written as
	\[
	\sup_{f\in\mathcal{F}_h}
	\big|(\mathbb{P}_N-P)f\big|
	=
	O_p\!\left(
	\sqrt{\frac{\log N}{N h^s}}
	\right).
	\]
	This reformulation is adopted throughout the remainder of the paper. If a preliminary or data-driven bandwidth exceeds one in finite samples, it may be truncated without affecting the function class or the asymptotic rates.
\end{remark}

\begin{proof}
	The argument follows standard empirical process techniques and is divided into three steps.

		\medskip
		\noindent\textbf{(i)}
		By (R8), there exists $C_\ell<\infty$ such that
		$|\ell(\psi^{-1}(\phi_p(\mathbf{Y}-\mathbf{y})^\top\boldsymbol{\gamma}),U)|\le C_\ell$ for all $(\mathbf{y},\boldsymbol{\gamma},U,\mathbf{Y})$.
		By (R2), $K$ is bounded and supported on $[-1,1]^s$, so $|K(v)|\le K_{\max}$ and
		$K((\mathbf{Y}-\mathbf{y})/h)\neq0$ implies $\|\mathbf{Y}-\mathbf{y}\|_\infty\le h$. Hence,
		$|f_{\mathbf{y},\boldsymbol{\gamma}}(U,\mathbf{Y})|\le C_\ell K_{\max} h^{-s}\mathbf 1\{\|\mathbf{Y}-\mathbf{y}\|_\infty\le h\}$.
		
		For each $\mathbf{y}\in\mathcal Y_0$, define the \emph{local} envelope
		$F_{\mathbf{y},h}(U,\mathbf{Y}):=C_\ell K_{\max} h^{-s}\mathbf 1\{\|\mathbf{Y}-\mathbf{y}\|_\infty\le h\}$.
		Then $|f_{\mathbf{y},\boldsymbol{\gamma}}(U,\mathbf{Y})|\le F_{\mathbf{y},h}(U,\mathbf{Y})$ for all $\boldsymbol{\gamma}\in\Gamma$.
		
		Moreover,
		\[
		\|F_{\mathbf{y},h}\|_{L_2(P)}^2
		=
		C_\ell^2 K_{\max}^2 h^{-2s}\,P(\|\mathbf{Y}-\mathbf{y}\|_\infty\le h).
		\]
		By (R1), $f_{\mathbf{Y}}$ is continuous and satisfies $c_f\le f_{\mathbf{Y}}\le C_f$ on an open neighborhood of
		$\mathcal Y_0$. Therefore, for all sufficiently small $h>0$ and uniformly over $\mathbf{y} \in\mathcal Y_0$,
		\[
		c_f\,\lambda_s([-h,h]^s)
		\le
		P(\|\mathbf{Y}-\mathbf{y}\|_\infty\le h)
		\le
		C_f\,\lambda_s([-h,h]^s),
		\]
		where $\lambda_s$ denotes Lebesgue measure on $\mathbb R^s$. Since $\lambda_s([-h,h]^s)=(2h)^s$,
		there exist constants $0<c_1\le c_2<\infty$ (independent of $h$) such that
		$c_1 h^s\le P(\|\mathbf{Y}-\mathbf{y}\|_\infty\le h)\le c_2 h^s$ uniformly for $\mathbf{y} \in\mathcal Y_0$.
		Consequently, there exist constants $0<C_1\le C_2<\infty$ (independent of $h$) such that
		\[
		C_1 h^{-s}\le \|F_{\mathbf{y},h}\|_{L_2(P)}^2\le C_2 h^{-s}
		\quad\text{uniformly for }\mathbf{y} \in\mathcal Y_0.
		\]
		Taking square roots yields
		\[
		\sqrt{C_1}\,h^{-s/2}\le \|F_{\mathbf{y},h}\|_{L_2(P)}\le \sqrt{C_2}\,h^{-s/2}
		\quad\text{uniformly for }\mathbf{y} \in\mathcal Y_0,
		\]
		which proves part \textnormal{(i)}.
		
		\medskip
		\noindent\textbf{(ii)}
		Define the normalized class
		\[
		\mathcal G_h
		:=
		\left\{
		g_{\mathbf{y},\boldsymbol{\gamma}}:=\frac{f_{\mathbf{y},\boldsymbol{\gamma}}}{\|F_{\mathbf{y},h}\|_{L_2(P)}}:\ \mathbf{y} \in\mathcal Y_0,\ \boldsymbol{\gamma}\in\Gamma
		\right\}.
		\]
		On the event $\|(\mathbf{Y}-\mathbf{y})/h\|_\infty\le1$, the vector $\phi_p(\mathbf{Y}-\mathbf{y})$ is uniformly bounded, so
		$\psi^{-1}(\phi_p(\mathbf{Y}-\mathbf{y})^\top\boldsymbol{\gamma})$ ranges over a compact subset of $\Theta_0$.
		Using a first-order Taylor expansion of
		$\ell(\psi^{-1}(\phi_p(\mathbf{Y}-\mathbf{y})^\top\boldsymbol{\gamma}),U)$ with respect to $(\mathbf{y},\boldsymbol{\gamma})$, together with
		the derivative bounds in (R4), (R5), and (R6), there exists $L_0>0$ such that
		\[
		|f_{\mathbf{y},\boldsymbol{\gamma}}(U,\mathbf{Y})-f_{\mathbf{y}',\boldsymbol{\gamma}'}(U,\mathbf{Y})|
		\le
		L_0 h^{-s}(\|\mathbf{Y}-\mathbf{y}'\|+\|\boldsymbol{\gamma}-\boldsymbol{\gamma}'\|)\,
		\mathbf 1\{\|\mathbf{Y}-\mathbf{y}\|_\infty\le h\ \text{or}\ \|\mathbf{Y}-\mathbf{y}'\|_\infty\le h\}.
		\]
		In particular, the indicator on the right-hand side is bounded by
		$\mathbf 1\{\|\mathbf{Y}-\mathbf{y}\|_\infty\le 2h\}$ whenever $\|\mathbf{Y}-\mathbf{y}'\|_\infty\le h$.
		Using (R1) as in part \textnormal{(i)} gives $P(\|\mathbf{Y}-\mathbf{y}\|_\infty\le 2h)\le C h^s$ uniformly in $\mathbf{y}$,
		and therefore there exists $C_1>0$ (independent of $h$) such that
		\[
		\|f_{\mathbf{y},\boldsymbol{\gamma}}-f_{\mathbf{y}',\boldsymbol{\gamma}'}\|_{L_2(P)}
		\le
		C_1 h^{-s/2}(\|\mathbf{Y}-\mathbf{y}'\|+\|\boldsymbol{\gamma}-\boldsymbol{\gamma}'\|).
		\]
		From part \textnormal{(i)}, $\|F_{\mathbf{y},h}\|_{L_2(P)}\ge \sqrt{C_1}\,h^{-s/2}$ uniformly in $\mathbf{y}$, hence
		\[
		\|g_{\mathbf{y},\boldsymbol{\gamma}}-g_{\mathbf{y}',\boldsymbol{\gamma}'}\|_{L_2(P)}
		=
		\frac{\|f_{\mathbf{y},\boldsymbol{\gamma}}-f_{\mathbf{y}',\boldsymbol{\gamma}'}\|_{L_2(P)}}{\|F_{\mathbf{y},h}\|_{L_2(P)}}
		\le
		C_2(\|\mathbf{Y}-\mathbf{y}'\|+\|\boldsymbol{\gamma}-\boldsymbol{\gamma}'\|),
		\]
		for a constant $C_2>0$ independent of $h$. Thus $\mathcal G_h$ is a finite-dimensional Lipschitz
		class with Lipschitz constant independent of $h$.
		
		Since $\mathcal Y_0\subset\mathbb R^s$ and $\Gamma\subset\mathbb R^{d_p}$ are compact,
		standard bounds on covering numbers for bounded subsets of Euclidean space imply that
		$\mathcal Y_0$ can be covered by at most $C_Y h^{-s}$ sup--norm balls of radius $h$ and
		$\Gamma$ by at most $C_\Gamma\varepsilon^{-d_p}$ Euclidean balls of radius $\varepsilon$
		(see Chapter~2.1--2.2 of \cite{van1996weak}). Combining these coverings yields
		\[
		N(\varepsilon,\mathcal G_h,L_2(P))\le C h^{-s}\varepsilon^{-d_p},
		\]
		for some constant $C>0$ independent of $h$.
		
		Finally, since $\|F_{\mathbf{y},h}\|_{L_2(P)}$ is bounded above and below by constant multiples of $h^{-s/2}$
		uniformly in $\mathbf{y}$, rescaling back to $\mathcal F_h$ yields
		\[
		\log N\!\Bigl(\varepsilon\,\|F_{\mathbf{y},h}\|_{L_2(P)},\ \mathcal F_h,\ L_2(P)\Bigr)
		\le
		C\bigl[1+\log(1/h)+\log(1/\varepsilon)\bigr],
		\]
		which proves part \textnormal{(ii)}.
		
		\medskip
		\noindent\textbf{(iii)}
		Let $\mathbb G_N=\sqrt N(\mathbb P_N-P)$. Applying Theorem~2.14.1 of \cite{van1996weak} to the
		$P$-measurable class $\mathcal F_h$ with envelope $F_{\mathbf{y},h}$ gives a constant $C'>0$ such that
		\[
		E\Big[\sup_{f\in\mathcal F_h}|\mathbb G_N f|\Big]
		\le
		C' J(\theta_N,\mathcal F_h)\,\sup_{\mathbf{y} \in\mathcal Y_0}\|F_{\mathbf{y},h}\|_{L_2(P)},
		\qquad
		\theta_N=\frac{\sup_{f\in\mathcal F_h}\|f\|_{L_2(P)}}{\sup_{\mathbf{y} \in\mathcal Y_0}\|F_{\mathbf{y},h}\|_{L_2(P)}}.
		\]
		By part \textnormal{(i)}, $\sup_{\mathbf{y} \in\mathcal Y_0}\|F_{\mathbf{y},h}\|_{L_2(P)}\le C h^{-s/2}$, and by the same
		localization argument as in part \textnormal{(i)}, $\sup_{f\in\mathcal F_h}\|f\|_{L_2(P)}\le C h^{-s/2}$,
		so $\theta_N\le C$ for all sufficiently small $h$.
		Using the entropy bound from part \textnormal{(ii)}, the entropy integral satisfies
		$J(\theta_N,\mathcal F_h)\le C\sqrt{\log(1/h)}$. Therefore,
		\[
		E\Big[\sup_{f\in\mathcal F_h}|\mathbb G_N f|\Big]
		\le
		C'' h^{-s/2}\sqrt{\log(1/h)}.
		\]
		It follows that
		\[
		\sup_{f\in\mathcal F_h}|(\mathbb P_N-P)f|
		=
		\frac{1}{\sqrt N}\sup_{f\in\mathcal F_h}|\mathbb G_N f|
		=
		O_p\!\left(\sqrt{\frac{\log(1/h)}{N h^s}}\right),
		\]
		which proves part \textnormal{(iii)}. Applying Remark~\ref{rem:logN-rate} we can replace $\log(1/h)$ with $\log N$.
	\end{proof}

	The following Corollary follows from Lemma~\ref{lem:kernel-class-entropy}.

\begin{corollary}\label{lem:uniform-LLN-rate}
	Assume that \textnormal{(R1)--(R8)} hold and that $\Gamma \subset \mathbb{R}^{d_p}$ is compact.
	Then, as $N \to \infty$,
	\begin{align}
		\sup_{\mathbf{y}\in\mathcal{Y}_0}\sup_{\boldsymbol{\gamma}\in\Gamma}
		\big|
		\tilde L_N(\boldsymbol{\gamma};\mathbf{y}) - Q(\boldsymbol{\gamma};\mathbf{y})
		\big|
		&=
		O_p\!\left(
		\sqrt{\frac{\log (1/h)}{N h^s}}
		\right),
		\label{eq:UNIF-OBJ-rate}
		\\[0.25cm]
		\sup_{\mathbf{y}\in\mathcal{Y}_0}\sup_{\boldsymbol{\gamma}\in\Gamma}
		\bigg\|
		\frac{1}{N h^s}\nabla_{\boldsymbol{\gamma}}L_N(\boldsymbol{\gamma};\mathbf{y})
		-
		E\!\left[
		\nabla_{\boldsymbol{\gamma}}\tilde L_N(\boldsymbol{\gamma};\mathbf{y})
		\right]
		\bigg\|
		&=
		O_p\!\left(
		\sqrt{\frac{\log (1/h)}{N h^s}}
		\right),
		\label{eq:UNIF-GRAD-rate}
		\\[0.25cm]
		\sup_{\mathbf{y}\in\mathcal{Y}_0}\sup_{\boldsymbol{\gamma}\in\Gamma}
		\bigg\|
		\frac{1}{N h^s}\nabla^2_{\boldsymbol{\gamma}} L_N(\boldsymbol{\gamma};\mathbf{y})
		-
		E\!\left[
		\nabla^2_{\boldsymbol{\gamma}}\tilde L_N(\boldsymbol{\gamma};\mathbf{y})
		\right]
		\bigg\|
		&=
		O_p\!\left(
		\sqrt{\frac{\log (1/h)}{N h^s}}
		\right).
		\label{eq:UNIF-HESS-rate}
	\end{align}
\end{corollary}

\noindent\emph{Remark.}
The stated rates follow from uniform laws of large numbers for kernel-weighted empirical processes.
The argument relies on the compact support and boundedness of the kernel $K$, the envelope conditions in
\textnormal{(R5)}, and the bandwidth requirement $N h^s / \log(1/h) \to \infty$.

	\begin{proof}
		Recall from Section~\ref{sec:setup} that the (kernel-weighted) local log-likelihood at $\mathbf{y}$ is
		\[
		L_N(\boldsymbol{\gamma};\mathbf{y})
		=
		\sum_{i=1}^N
		\ell\!\Bigl(
		\psi^{-1}\bigl(\phi_p(Y_i-\mathbf{y})^\top\boldsymbol{\gamma}\bigr),\,U_i
		\Bigr)\,
		K\!\left(\frac{Y_i-\mathbf{y}}{h}\right),
		\]
		and define the scaled objective and its expectation by
		$\tilde L_N(\boldsymbol{\gamma};\mathbf{y})=\frac{1}{N h^s} L_N(\boldsymbol{\gamma};\mathbf{y})$ and
		$Q(\boldsymbol{\gamma};\mathbf{y})=E[\tilde L_N(\boldsymbol{\gamma};\mathbf{y})]$.
		Let $\mathbb{P}_N f = N^{-1}\sum_{i=1}^N f(U_i,Y_i)$ and $P f = E[f(U,\mathbf{Y})]$ denote
		the empirical and population measures. For $(\mathbf{y},\boldsymbol{\gamma})\in\mathcal{Y}_0\times\Gamma$ define
		\[
		f_{\mathbf{y},\boldsymbol{\gamma}}(U,\mathbf{Y})
		:=
		h^{-s}\,
		\ell\!\Bigl(
		\psi^{-1}\bigl(\phi_p(\mathbf{Y}-\mathbf{y})^\top\boldsymbol{\gamma}\bigr),\,U
		\Bigr)\,
		K\!\left(\frac{\mathbf{Y}-\mathbf{y}}{h}\right).
		\]
		Then $\tilde L_N(\boldsymbol{\gamma};\mathbf{y})=\mathbb{P}_N f_{\mathbf{y},\boldsymbol{\gamma}}$ and
		$Q(\boldsymbol{\gamma};\mathbf{y})=P f_{\mathbf{y},\boldsymbol{\gamma}}$, so
		\[
		\sup_{\mathbf{y}\in\mathcal{Y}_0}\sup_{\boldsymbol{\gamma}\in\Gamma}
		\bigl|\tilde L_N(\boldsymbol{\gamma};\mathbf{y})-Q(\boldsymbol{\gamma};\mathbf{y})\bigr|
		=
		\sup_{f\in\mathcal{F}_h}\bigl|(\mathbb{P}_N-P)f\bigr|,
		\qquad
		\mathcal{F}_h=\{f_{\mathbf{y},\boldsymbol{\gamma}}:\mathbf{y}\in\mathcal{Y}_0,\ \boldsymbol{\gamma}\in\Gamma\}.
		\]
		The uniform rate \eqref{eq:UNIF-OBJ-rate} therefore follows directly from
		Lemma~\ref{lem:kernel-class-entropy}(iii).
		
		\medskip\noindent
		\textit{Gradient bound.}
		By the chain rule (see Section~\ref{sec:setup}),
		\[
		\nabla_{\boldsymbol{\gamma}}L_N(\boldsymbol{\gamma};\mathbf{y})
		=
		\sum_{i=1}^N
		\ell'\!\Bigl(
		\psi^{-1}(\phi_p(Y_i-\mathbf{y})^\top\boldsymbol{\gamma}),\,U_i
		\Bigr)\,
		(\psi^{-1})'\!\bigl(\phi_p(Y_i-\mathbf{y})^\top\boldsymbol{\gamma}\bigr)\,
		\phi_p(Y_i-\mathbf{y})\,
		K\!\left(\frac{Y_i-\mathbf{y}}{h}\right),
		\]
		and hence
		\[
		\frac{1}{N h^s}\,\nabla_{\boldsymbol{\gamma}}L_N(\boldsymbol{\gamma};\mathbf{y})
		=
		\mathbb{P}_N g_{\mathbf{y},\boldsymbol{\gamma}},
		\]
		where each coordinate of $g_{\mathbf{y},\boldsymbol{\gamma}}(U,\mathbf{Y})$ is
		\begin{equation}\label{eq:gradclass-def}
			g_{\mathbf{y},\boldsymbol{\gamma}}(U,\mathbf{Y})
			=
			h^{-s}\,
			\ell'\!\Bigl(
			\psi^{-1}(\phi_p(\mathbf{Y}-\mathbf{y})^\top\boldsymbol{\gamma}),\,U
			\Bigr)\,
			(\psi^{-1})'\!\bigl(\phi_p(\mathbf{Y}-\mathbf{y})^\top\boldsymbol{\gamma}\bigr)\,
			\phi_p(\mathbf{Y}-\mathbf{y})\,
			K\!\left(\frac{\mathbf{Y}-\mathbf{y}}{h}\right).
		\end{equation}
		Let $\mathcal{G}_h := \{g_{\mathbf{y},\boldsymbol{\gamma}}:\mathbf{y}\in\mathcal{Y}_0,\ \boldsymbol{\gamma}\in\Gamma\}$.
		
		By (R4)--(R6) and the envelope condition in (R5), there exists $C_0>0$ such that
		$|\ell'(\theta,U)|\le C_0$ on the relevant range (in particular, $E[M_1(U)^2]<\infty$
		suffices for the $L_2$ arguments below). By (R2), $K((\mathbf{Y}-\mathbf{y})/h)\neq0$ implies
		$\|\mathbf{Y}-\mathbf{y}\|_\infty\le h$, hence $\|\mathbf{Y}-\mathbf{y}\|\le C_K h$ for some $C_K>0$, and therefore
		\[
		|g_{\mathbf{y},\boldsymbol{\gamma}}(U,\mathbf{Y})|
		\le
		C_0 h^{-s}\mathbf{1}\{\|\mathbf{Y}-\mathbf{y}\|\le C_K h\}.
		\]
		
		Define the (global) envelope
		\[
		F_h^{(1)}(Y)
		:=
		C_0 h^{-s}\mathbf{1}\{\operatorname{dist}_\infty(Y,\mathcal{Y}_0)\le C_K h\}.
		\]
		Since $\mathcal{Y}_0$ is compact and has nonempty interior and $f_{\mathbf{Y}}$ is bounded away from
		zero and above on a neighborhood of $\mathcal{Y}_0$ by (R1), there exist constants
		$0<c_1\le c_2<\infty$ such that for all sufficiently small $h$,
		\[
		c_1
		\le
		P\!\bigl(\operatorname{dist}_\infty(Y,\mathcal{Y}_0)\le C_K h\bigr)
		\le
		c_2,
		\]
		and consequently
		\[
		C_1 h^{-s}
		\le
		\|F_h^{(1)}\|_{L_2(P)}
		\le
		C_2 h^{-s}
		\]
		for some constants $0<C_1\le C_2<\infty$.
		
		For later use, note that for each fixed $\mathbf{y}\in\mathcal{Y}_0$,
		\[
		\|g_{\mathbf{y},\boldsymbol{\gamma}}\|_{L_2(P)}^2
		\le
		C_0^2 h^{-2s} P(\|\mathbf{Y}-\mathbf{y}\|\le C_K h),
		\]
		and by (R1) there exist constants $0<\underline c\le \overline c<\infty$ such that
		$\underline c\,h^s \le P(\|\mathbf{Y}-\mathbf{y}\|\le C_K h)\le \overline c\,h^s$ uniformly in $\mathbf{y}\in\mathcal{Y}_0$.
		Hence there exists $C_3>0$ such that
		\[
		\sup_{\mathbf{y}\in\mathcal{Y}_0}\sup_{\boldsymbol{\gamma}\in\Gamma}\|g_{\mathbf{y},\boldsymbol{\gamma}}\|_{L_2(P)} \le C_3 h^{-s/2}.
		\]
		
		A first-order Taylor expansion of \eqref{eq:gradclass-def} in $(\mathbf{y},\boldsymbol{\gamma})$, together with
		(R4)--(R6), yields for some $L_0>0$,
		\(|g_{\mathbf{y},\boldsymbol{\gamma}}(U,\mathbf{Y})-g_{\mathbf{y}',\boldsymbol{\gamma}'}(U,\mathbf{Y})|
		\le
		L_0 h^{-s}(\|\mathbf{Y}-\mathbf{y}'\|+\|\boldsymbol{\gamma}-\boldsymbol{\gamma}'\|)
		\mathbf{1}\{\operatorname{dist}_\infty(Y,\mathcal{Y}_0)\le 2C_K h\}.
		\)
		Since $P(\operatorname{dist}_\infty(Y,\mathcal{Y}_0)\le 2C_K h)$ is bounded away from $0$
		and $1$ for small $h$, this implies
		\(\|g_{\mathbf{y},\boldsymbol{\gamma}}-g_{\mathbf{y}',\boldsymbol{\gamma}'}\|_{L_2(P)}
		\le
		C h^{-s}(\|\mathbf{Y}-\mathbf{y}'\|+\|\boldsymbol{\gamma}-\boldsymbol{\gamma}'\|)
		\)
		for a constant $C>0$. Combining this with $\|F_h^{(1)}\|_{L_2(P)}\ge C_1 h^{-s}$ shows that
		the normalized class is Lipschitz in $(\mathbf{y},\boldsymbol{\gamma})$ with a constant independent of $h$.
		
		Cover $\mathcal{Y}_0$ by a sup-norm grid of mesh $h$ (at most $C_Y h^{-s}$ points)
		and $\Gamma$ by a Euclidean $\delta$-net (at most $C_\Gamma\delta^{-d_p}$ points).
		The Lipschitz bound yields
		\[
		N(\delta,\mathcal{G}_h,L_2(P))
		\le
		C h^{-s}\delta^{-d_p}.
		\]
		
		Now apply a maximal inequality for empirical processes with bounded envelopes and
		polynomial entropy (e.g. Theorem~2.14.1 in \cite{van1996weak}) to $\mathcal{G}_h$.
		With $\theta_N=\sup_{g\in\mathcal{G}_h}\|g\|_{L_2(P)}/\|F_h^{(1)}\|_{L_2(P)}\le (C_3/C_1)h^{s/2}$
		and the entropy bound above, we obtain $J(\theta_N,\mathcal{G}_h)\le C\,\theta_N\sqrt{\log(1/h)}$,
		and hence
		\[
		E\Big[\sup_{g\in\mathcal{G}_h}|\mathbb{G}_N g|\Big]
		\le
		C\,J(\theta_N,\mathcal{G}_h)\,\|F_h^{(1)}\|_{L_2(P)}
		\le
		C' h^{-s/2}\sqrt{\log(1/h)}.
		\]
		Therefore,
		\[
		\sup_{g\in\mathcal{G}_h}|(\mathbb{P}_N-P)g|
		=
		O_p\!\left(\sqrt{\frac{\log(1/h)}{N h^s}}\right).
		\]
		Since
		$\frac{1}{N h^s}\nabla_{\boldsymbol{\gamma}}L_N(\boldsymbol{\gamma};\mathbf{y})-E[\nabla_{\boldsymbol{\gamma}}\tilde L_N(\boldsymbol{\gamma};\mathbf{y})]
		=(\mathbb{P}_N-P)g_{\mathbf{y},\boldsymbol{\gamma}}$,
		taking suprema yields \eqref{eq:UNIF-GRAD-rate}.
		
		\medskip\noindent
		\textit{Hessian bound.}
		Differentiating once more,
		\[
		\nabla^2_{\boldsymbol{\gamma}} L_N(\boldsymbol{\gamma};\mathbf{y})
		=
		\sum_{i=1}^N H(U_i,Y_i;\boldsymbol{\gamma},\mathbf{y})\,
		K\!\left(\frac{Y_i-\mathbf{y}}{h}\right),
		\]
		where each entry of $H(U,Y;\boldsymbol{\gamma},\mathbf{y})$ is a linear combination of terms involving
		$\ell''$, $\ell'$, $(\psi^{-1})'$, $(\psi^{-1})''$, and $\phi_p(\mathbf{Y}-\mathbf{y})\phi_p(\mathbf{Y}-\mathbf{y})^\top$.
		Using (R4)--(R6) and the envelope bounds in (R5) for $\ell'$ and $\ell''$, each scalar entry
		$h_{\mathbf{y},\boldsymbol{\gamma}}(U,\mathbf{Y})$ of $(Nh^s)^{-1}\nabla^2_{\boldsymbol{\gamma}} L_N(\boldsymbol{\gamma};\mathbf{y})=\mathbb{P}_N h_{\mathbf{y},\boldsymbol{\gamma}}$
		satisfies $|h_{\mathbf{y},\boldsymbol{\gamma}}(U,\mathbf{Y})|\le C h^{-s}\mathbf{1}\{\|\mathbf{Y}-\mathbf{y}\|\le C_K h\}$.
		Let $\mathcal{H}_h:=\{h_{\mathbf{y},\boldsymbol{\gamma}}:\mathbf{y}\in\mathcal{Y}_0,\ \boldsymbol{\gamma}\in\Gamma\}$.
		The same covering and maximal inequality argument as for $\mathcal{G}_h$
		(with $\sup_{h\in\mathcal{H}_h}\|h\|_{L_2(P)}\le C h^{-s/2}$) yields
		\[
		\sup_{h\in\mathcal{H}_h}|(\mathbb{P}_N-P)h|
		=
		O_p\!\left(\sqrt{\frac{\log(1/h)}{N h^s}}\right),
		\]
		which implies \eqref{eq:UNIF-HESS-rate}.
		
		Combining the three bounds completes the proof. In all the three cases, by applying Remark~\ref{rem:logN-rate} we can replace $\log(1/h)$ with $\log N$.
	\end{proof}

\begin{lemma}\label{lem:score-expansion}
	Assume that \textnormal{(R1)--(R7)} hold and that the local polynomial order $p$ is odd.
	Let $\boldsymbol{\gamma}^\ast(\mathbf{y})$ denote the unique maximizer of the population
	criterion $Q(\cdot;\mathbf{y})$, and define
	\[
	S_N(\mathbf{y})
	=
	\frac{1}{N h^s}\,
	\nabla_{\boldsymbol{\gamma}}L_N(\boldsymbol{\gamma}^\ast(\mathbf{y});\mathbf{y}).
	\]
	Then, as $N \to \infty$,
	\[
	\sup_{\mathbf{y}\in\mathcal{Y}_0}
	\big\| S_N(\mathbf{y}) \big\|
	=
	O_p\!\left(
	h^{p+1}
	+
	\sqrt{\frac{\log (1/h)}{N h^s}}
	\right).
	\]
	Under the additional condition stated in Remark~\ref{rem:logN-rate}, the logarithmic term
	$\log(1/h)$ may be replaced by $\log N$ for sufficiently large $N$.
\end{lemma}

	\begin{proof}
		We decompose the local score into its population part (bias) and empirical fluctuation, then bound each term uniformly in $\mathbf{y}$.
		
		Recall \(L_N(\boldsymbol{\gamma};\mathbf{y})
		=
		\sum_{i=1}^N
		\ell\!\Bigl(
		\psi^{-1}\big(\phi_p(Y_i - y)^\top \boldsymbol{\gamma}\big),\,U_i
		\Bigr)\,
		K\!\left(\frac{Y_i - y}{h}\right),
		\)
		and \(\tilde L_N(\boldsymbol{\gamma};\mathbf{y})
		=
		\frac{1}{N h^s} L_N(\boldsymbol{\gamma};\mathbf{y}),
		\\
		Q(\boldsymbol{\gamma};\mathbf{y})
		=
		E\big[\tilde L_N(\boldsymbol{\gamma};\mathbf{y})\big].
		\) Define the \emph{population score} at $\boldsymbol{\gamma}^\ast(\mathbf{y})$ by
		\[
		s_h(\mathbf{y})
		:=
		E\big[S_N(\mathbf{y})\big]
		=
		E\!\left[\frac{1}{N h^s}\,\nabla_{\boldsymbol{\gamma}}L_N(\boldsymbol{\gamma}^\ast(\mathbf{y});\mathbf{y})\right]
		=
		E\big[\nabla_{\boldsymbol{\gamma}}\tilde L_N(\boldsymbol{\gamma}^\ast(\mathbf{y});\mathbf{y})\big].
		\]
		Then we can write \(S_N(\mathbf{y})
		=
		\underbrace{\bigl(S_N(\mathbf{y}) - s_h(\mathbf{y})\bigr)}_{\text{empirical fluctuation}}
		+
		\underbrace{s_h(\mathbf{y})}_{\text{population term}}.
		\) Hence
		\[
		\sup_{\mathbf{y}\in\mathcal{Y}_0} \|S_N(\mathbf{y})\|
		\;\le\;
		\sup_{\mathbf{y}\in\mathcal{Y}_0} \|S_N(\mathbf{y})-s_h(\mathbf{y})\|
		+
		\sup_{\mathbf{y}\in\mathcal{Y}_0} \|s_h(\mathbf{y})\|.
		\]
		We bound these two terms separately.
		
		By definition,
		\[
		S_N(\mathbf{y})
		=
		\frac{1}{N h^s}\,\nabla_{\boldsymbol{\gamma}}L_N(\boldsymbol{\gamma}^\ast(\mathbf{y});\mathbf{y})
		=
		\mathbb{P}_N g_{\mathbf{y},\boldsymbol{\gamma}^\ast(\mathbf{y})},
		\]
		where, componentwise, $g_{\mathbf{y},\boldsymbol{\gamma}}(U,\mathbf{Y})$ is as given by \eqref{eq:gradclass-def} and $\mathbb{P}_N$ is the empirical measure based on $(U_i,Y_i)_{i=1}^N$. Let \(\mathcal{G}_h
		=
		\{ g_{\mathbf{y},\boldsymbol{\gamma}} : \mathbf{y}\in\mathcal{Y}_0,\ \boldsymbol{\gamma}\in\Gamma\}.
		\)
		By the same envelope, Lipschitz, and entropy arguments used in Lemma~\ref{lem:kernel-class-entropy} (see the proof of Corollary~\ref{lem:uniform-LLN-rate}), the class $\mathcal{G}_h$ satisfies
		\[
		\sup_{\mathbf{y}\in\mathcal{Y}_0}\sup_{\boldsymbol{\gamma}\in\Gamma}
		\big\|
		\mathbb{P}_N g_{\mathbf{y},\boldsymbol{\gamma}} - P g_{\mathbf{y},\boldsymbol{\gamma}}
		\big\|
		=
		O_p\!\left(
		\sqrt{\frac{\log(1/h)}{N h^s}}
		\right).
		\]
		In particular, since $\{\boldsymbol{\gamma}^\ast(\mathbf{y}): \mathbf{y}\in\mathcal{Y}_0\}\subset\Gamma$, we have
		\begin{align}
			\sup_{\mathbf{y}\in\mathcal{Y}_0}
			\big\|S_N(\mathbf{y}) - s_h(\mathbf{y})\big\|
			=&
			\sup_{\mathbf{y}\in\mathcal{Y}_0}
			\big\|
			\mathbb{P}_N g_{\mathbf{y},\boldsymbol{\gamma}^\ast(\mathbf{y})} - P g_{\mathbf{y},\boldsymbol{\gamma}^\ast(\mathbf{y})}
			\big\| \nonumber\\
			\le&
			\sup_{\mathbf{y}\in\mathcal{Y}_0}\sup_{\boldsymbol{\gamma}\in\Gamma}
			\big\|
			\mathbb{P}_N g_{\mathbf{y},\boldsymbol{\gamma}} - P g_{\mathbf{y},\boldsymbol{\gamma}}
			\big\|
			=
			O_p\!\left(
			\sqrt{\frac{\log(1/h)}{N h^s}}
			\right).  \label{bound_empirical_fluctuation} 
		\end{align}

		We now show that $\sup_{\mathbf{y}\in\mathcal{Y}_0} \|s_h(\mathbf{y})\| = O(h^{p+1})$ uniformly in $\mathbf{y}$ with $p$ odd. Write $s_h(\mathbf{y})$ more explicitly:
		\[
		s_h(\mathbf{y})
		=
		E\big[\nabla_{\boldsymbol{\gamma}}\tilde L_N(\boldsymbol{\gamma}^\ast(\mathbf{y});\mathbf{y})\big]
		=
		E\!\left[
		h^{-s}
		\nabla_{\boldsymbol{\gamma}}\ell\!\Bigl(
		\psi^{-1}(\phi_p(Y - \mathbf{y})^\top \boldsymbol{\gamma}^\ast(\mathbf{y})),\,U
		\Bigr)\,
		K\!\left(\frac{Y - \mathbf{y}}{h}\right)
		\right].
		\]
		
		Conditioning on $\mathbf{Y}$ and using the chain rule, each coordinate of the integrand is a bounded linear combination (by (R4)--(R6)) of terms of the form
		\[
		h^{-s}
		E\!\left[
		\ell'\!\Bigl(\theta_{\mathbf{y}}^\ast(Y),U\Bigr) \mid Y
		\right]
		(\psi^{-1})'\big(\phi_p(\mathbf{Y}-\mathbf{y})^\top\boldsymbol{\gamma}^\ast(\mathbf{y})\big)
		\phi_p(\mathbf{Y}-\mathbf{y})
		K\!\left(\frac{\mathbf{Y}-\mathbf{y}}{h}\right),
		\]
		where \(\theta_{\mathbf{y}}^\ast(Y)
		=
		\psi^{-1}\big(\phi_p(\mathbf{Y}-\mathbf{y})^\top\boldsymbol{\gamma}^\ast(\mathbf{y})\big).
		\)  Let $ \theta(Y) = \psi^{-1}(\nu(Y))$ denote the true calibration evaluated at $\mathbf{Y}$, so that $U\mid Y$ has copula parameter $\theta(Y)$ and, in particular,
		\(
		E\big[\,\ell'(\theta(Y),U)\mid Y\big] = 0
		\) since the conditional copula of $U \mid Y=\mathbf{y}$ is $C(\,\cdot \mid \theta(\mathbf{y})).$
		Hence we may write
		\(
		E\big[\,\ell'(\theta_{\mathbf{y}}^\ast(Y),U)\mid Y\big]
		=
		E\big[\,\ell'(\theta_{\mathbf{y}}^\ast(Y),U) - \ell'(\theta(Y),U)\mid Y\big].
		\)
		
		By (R5) and (R6), $\ell'$ is differentiable in $\theta$ with uniformly bounded derivative, and $\psi^{-1}$ has uniformly bounded derivatives on the relevant range. Thus,
		\[
		\big|
		E\big[\,\ell'(\theta_{\mathbf{y}}^\ast(Y),U)\mid Y\big]
		\big|
		\;\le\;
		C\,|\theta_{\mathbf{y}}^\ast(Y) - \theta(Y)|
		\]
		for some constant $C<\infty$. Furthermore, by definition of $\boldsymbol{\gamma}^\ast(\mathbf{y})$ and (R4), $\phi_p(\cdot)$ is the order-$p$ local polynomial basis and $\nu$ is $(p+2)$-times continuously differentiable. A standard local polynomial approximation argument (see, e.g., \cite{fan1996local}) then yields the uniform expansion
		\[
		\nu(Y)
		=
		\phi_p(\mathbf{Y}-\mathbf{y})^\top \boldsymbol{\gamma}^\ast(\mathbf{y})
		+
		R_{p+1}(\mathbf{Y},\mathbf{y}),
		\]
		with remainder satisfying
		\[
		\sup_{\mathbf{y}\in\mathcal{Y}_0}
		\sup_{\|\mathbf{Y}-\mathbf{y}\|_\infty\le Ch}
		\big|R_{p+1}(Y,\mathbf{y})\big|
		\le
		C' \|\mathbf{Y}-\mathbf{y}\|^{p+1}
		\]
		for some $C',C<\infty$. Since $\psi^{-1}$ has bounded first derivative (R6), this implies
		\[
		|\theta_{\mathbf{y}}^\ast(Y) - \theta(Y)|
		=
		\big|\psi^{-1}(\phi_p(\mathbf{Y}-\mathbf{y})^\top \boldsymbol{\gamma}^\ast(\mathbf{y})) - \psi^{-1}(\nu(Y))\big|
		\;\le\;
		C'' \|\mathbf{Y}-\mathbf{y}\|^{p+1}
		\]
		for $\|\mathbf{Y}-\mathbf{y}\|_\infty\le Ch$ and some $C''<\infty$. Putting these bounds together and using boundedness of $(\psi^{-1})'$, $\phi_p(\cdot)$ and $K$ on the local region $\|\mathbf{Y}-\mathbf{y}\|_\infty\le Ch$, we obtain (for each coordinate)
		\[
		\big\|s_h(\mathbf{y})\big\|
		\;\le\;
		C_1
		E\!\left[
		h^{-s}\|\mathbf{Y}-\mathbf{y}\|^{p+1}
		\mathbf{1}\{\|\mathbf{Y}-\mathbf{y}\|_\infty \le C h\}
		\right].
		\]
		
		Using (R1) and a change of variables $v = (\mathbf{Y}-\mathbf{y})/h$, we get
		\[
		\begin{aligned}
			\|s_h(\mathbf{y})\|
			&\le
			C_2
			\int_{\|z-y\|_\infty \le C h}
			h^{-s}\,\|z-y\|^{p+1}\, f_{\mathbf{Y}}(z)\,dz \\
			&=
			C_2
			\int_{\|v\|_\infty \le C}
			h^{-s}\,\|h v\|^{p+1}\, f_{\mathbf{Y}}(y+h v)\, h^s\, dv \\
			&=
			C_3 h^{p+1}
			\int_{\|v\|_\infty \le C}
			\|v\|^{p+1}\, dv \;\le\;
			C_4 h^{p+1}.
		\end{aligned}
		\]
		where the constants $C_j$ do not depend on $\mathbf{y}$ (by boundedness and continuity of
		$f_{\mathbf{Y}}$ on $\mathcal{Y}_0$ and of the derivatives in
		\textnormal{(R4)--(R6)}). Since $\mathcal{Y}_0$ is compact, this bound is uniform in $\mathbf{y}$, namely, \begin{align}\label{bound_population_term}
			\sup_{\mathbf{y}\in\mathcal{Y}_0} \, \|s_h(\mathbf{y})\|\;\le\;
			C_4 h^{p+1}.
		\end{align}
		
		The assumption that $p$ is odd and $K$ is symmetric ensures that all polynomial terms up to order $p$ vanish after integration (because the corresponding odd moments of $K$ are zero), so the first nonzero contribution is of order $h^{p+1}$, as above.

		Combining the bounds \eqref{bound_empirical_fluctuation} and~\eqref{bound_population_term}, we obtain
		\[
		\sup_{\mathbf{y}\in\mathcal{Y}_0} \|S_N(\mathbf{y})\|
		\;\le\;
		\sup_{\mathbf{y}\in\mathcal{Y}_0} \|S_N(\mathbf{y})-s_h(\mathbf{y})\|
		+
		\sup_{\mathbf{y}\in\mathcal{Y}_0} \|s_h(\mathbf{y})\|
		=
		O_p\!\left(
		\sqrt{\frac{\log(1/h)}{N h^s}}
		\right)
		+
		O\!\big(h^{p+1}\big),
		\]
		which proves the claim.
	\end{proof}

\begin{lemma}\label{lem:gamma-uniform-rate}
	Let $\hat{\boldsymbol{\gamma}}(\mathbf{y})$ be a local maximizer of
	$L_N(\boldsymbol{\gamma};\mathbf{y})$ over the compact set $\Gamma$, and let
	$\boldsymbol{\gamma}^\ast(\mathbf{y})$ denote the unique maximizer of the population
	criterion $Q(\cdot;\mathbf{y})$ guaranteed by \textnormal{(R7)}.
	Assume that \textnormal{(R1)--(R7)} hold and that the polynomial order $p$ is odd.
	Then
	\[
	\sup_{\mathbf{y}\in\mathcal{Y}_0}
	\big\|
	\hat{\boldsymbol{\gamma}}(\mathbf{y}) - \boldsymbol{\gamma}^\ast(\mathbf{y})
	\big\|
	=
	O_p\!\left(
	h^{p+1}
	+
	\sqrt{\frac{\log (1/h)}{N h^s}}
	\right).
	\]
	Under the additional condition in Remark~\ref{rem:logN-rate}, the logarithmic term
	$\log(1/h)$ may be replaced by $\log N$ for all sufficiently large $N$.
\end{lemma}

	\begin{proof}
		The proof is based on a Taylor expansion of the empirical score around 
		$\boldsymbol{\gamma}^\ast(\mathbf{y})$, together with the uniform stochastic expansion of the score 
		and uniform control of the Hessian.
		
		For each $\mathbf{y}\in\mathcal{Y}_0$, the local maximizer $\hat{\boldsymbol{\gamma}}(\mathbf{y})$ of $L_N(\cdot;\mathbf{y})$ 
		satisfies the score equation
		\(\nabla_{\boldsymbol{\gamma}}L_N(\hat{\boldsymbol{\gamma}}(\mathbf{y});\mathbf{y}) = 0.
		\)
		Recall the scaled score at $\gamma^\ast(\mathbf{y})$,
		\(S_N(\mathbf{y})
		=
		\frac{1}{N h^s}\,\nabla_\gamma L_N(\gamma^\ast(\mathbf{y});\mathbf{y}),
		\)
		as in Lemma~\ref{lem:score-expansion}, and define the scaled empirical Hessian by
		\(H_N(\gamma;\mathbf{y})
		:=
		\frac{1}{N h^s}\,\nabla_\gamma^2 L_N(\gamma;\mathbf{y}),
		\
		H(\gamma;\mathbf{y})
		:=
		E\big[H_N(\gamma;\mathbf{y})\big]
		=
		E\big[\nabla_\gamma^2 \tilde L_N(\gamma;\mathbf{y})\big].
		\) Fix $\mathbf{y}\in\mathcal{Y}_0$ and apply a mean-value expansion of 
		$\nabla_\gamma L_N(\gamma;\mathbf{y})$ (equivalently, of its scaled version) around 
		$\gamma^\ast(\mathbf{y})$:
		\[
		\begin{aligned}
			0
			=
			\frac{1}{N h^s}\,\nabla_\gamma L_N(\hat\gamma(\mathbf{y});\mathbf{y})
			=&
			\frac{1}{N h^s}\,\nabla_\gamma L_N(\gamma^\ast(\mathbf{y});\mathbf{y})
			+\\ &\qquad
			\Bigg[
			\int_0^1 H_N\big(\gamma^\ast(\mathbf{y}) + t(\hat\gamma(\mathbf{y})-\gamma^\ast(\mathbf{y}));\mathbf{y}\big)\,dt
			\Bigg]
			\big(\hat\gamma(\mathbf{y})-\gamma^\ast(\mathbf{y})\big).    
		\end{aligned}
		\]
		In other words, \(0
		=
		S_N(\mathbf{y})
		+
		H_N^\ast(\mathbf{y})\big(\hat\gamma(\mathbf{y})-\gamma^\ast(\mathbf{y})\big),
		\)
		where we define the average Hessian along the line segment between 
		$\gamma^\ast(\mathbf{y})$ and $\hat\gamma(\mathbf{y})$ by
		\[
		H_N^\ast(\mathbf{y})
		:=
		\int_0^1 H_N\big(\gamma^\ast(\mathbf{y}) + t(\hat\gamma(\mathbf{y})-\gamma^\ast(\mathbf{y}));\mathbf{y}\big)\,dt.
		\]
		Rearranging gives
		\begin{equation}
			\label{eq:gamma-diff-Taylor}
			\hat\gamma(\mathbf{y}) - \gamma^\ast(\mathbf{y})
			=
			-\big[H_N^\ast(\mathbf{y})\big]^{-1} S_N(\mathbf{y}),
		\end{equation}
		provided $H_N^\ast(\mathbf{y})$ is invertible.
		
		Thus, to bound $\|\hat\gamma(\mathbf{y})-\gamma^\ast(\mathbf{y})\|$ uniformly in $\mathbf{y}$, 
		it suffices to obtain a uniform bound on $\|S_N(\mathbf{y})\|$ and a uniform lower bound 
		on the eigenvalues of $-H_N^\ast(\mathbf{y})$.

		By Lemma~\ref{lem:score-expansion}, under \textnormal{(R1)--(R7)} with $p$ odd,
		\[
		\sup_{\mathbf{y}\in\mathcal{Y}_0} \big\| S_N(\mathbf{y}) \big\|
		=
		O_p\!\left( h^{p+1} + \sqrt{\frac{\log (1/h)}{N h^s}} \right).
		\]
		
		To control the Hessian uniformly, Corollary~\ref{lem:uniform-LLN-rate} yields
		\[
		\sup_{\mathbf{y}\in\mathcal{Y}_0}\sup_{\gamma\in\Gamma}
		\big\| H_N(\gamma;\mathbf{y}) - H(\gamma;\mathbf{y}) \big\|
		=
		O_p\!\left( \sqrt{\frac{\log (1/h)}{N h^s}} \right).
		\]
		
		On the other hand, Assumption~(R7) implies that the population Hessian is
		uniformly negative definite:
		\[
		v^\top [-H(\gamma;\mathbf{y})]\,v \;\ge\; c_H \|v\|^2,
		\qquad
		\text{for all } v\in\mathbb{R}^{d_p},\ (\gamma,\mathbf{y})\in\Gamma\times\mathcal{Y}_0.
		\]
		
		Since 
		\(
		\sqrt{\log(1/h)/(N h^s)} \to 0
		\)
		under (R3), it follows that
		\[
		\sup_{\mathbf{y}\in\mathcal{Y}_0}\sup_{\gamma\in\Gamma}
		\big\| H_N(\gamma;\mathbf{y}) - H(\gamma;\mathbf{y}) \big\|
		=o_p(1).
		\]
		Therefore, for any $\eta \in (0,c_H)$, there exists $N_0$ such that, for all 
		$N\ge N_0$,
		\[
		\mathbb{P}
		\!\left(
		\sup_{\mathbf{y}\in\mathcal{Y}_0}\sup_{\gamma\in\Gamma}
		\big\| H_N(\gamma;\mathbf{y}) - H(\gamma;\mathbf{y}) \big\| \le \eta
		\right)
		\to 1.
		\]
		
		On this event, Weyl's inequality \cite{bhatia1997matrix} (or an elementary perturbation bound) implies that
		\[
		v^\top[-H_N(\gamma;\mathbf{y})]\,v
		\;\ge\;
		(c_H - \eta)\,\|v\|^2,
		\qquad
		\forall v,\ (\gamma,\mathbf{y})\in\Gamma\times\mathcal{Y}_0.
		\]
		Choosing $\eta = c_H/2$, we obtain, with probability tending to one,
		\[
		v^\top[-H_N(\gamma;\mathbf{y})]\,v
		\;\ge\;
		\frac{c_H}{2}\,\|v\|^2,
		\qquad
		\forall v,\ (\gamma,\mathbf{y})\in\Gamma\times\mathcal{Y}_0.
		\]
		
		Since $\gamma^\ast(\mathbf{y})$ and $\hat\gamma(\mathbf{y})$ both lie in $\Gamma$, every point
		on the segment $\gamma^\ast(\mathbf{y}) + t(\hat\gamma(\mathbf{y})-\gamma^\ast(\mathbf{y}))$ for $t\in[0,1]$
		also lies in $\Gamma$ (because $\Gamma$ is convex; if needed, one may take $\Gamma$ 
		as the closed convex hull of the original compact set without affecting (R7)). 
		Therefore, on the same high-probability event,
		\[
		v^\top[-H_N^\ast(\mathbf{y})]\,v
		=
		\int_0^1 v^\top[-H_N(\gamma^\ast(\mathbf{y})+t(\hat\gamma(\mathbf{y})-\gamma^\ast(\mathbf{y}));\mathbf{y})]\,v\,dt
		\;\ge\;
		\frac{c_H}{2}\,\|v\|^2,
		\qquad \forall v,\ \mathbf{y}\in\mathcal{Y}_0.
		\]
		This shows that $H_N^\ast(\mathbf{y})$ is invertible for all $\mathbf{y}\in\mathcal{Y}_0$ and that
		\[
		\sup_{\mathbf{y}\in\mathcal{Y}_0}
		\big\|\big[H_N^\ast(\mathbf{y})\big]^{-1}\big\|
		\;\le\;
		\frac{2}{c_H}
		\quad\text{with probability tending to one.}
		\]

		Returning to the expansion \eqref{eq:gamma-diff-Taylor}, on the event where
		$H_N^\ast(\mathbf{y})$ is invertible for all $\mathbf{y}$, we have
		\(\hat\gamma(\mathbf{y}) - \gamma^\ast(\mathbf{y})
		=
		-\big[H_N^\ast(\mathbf{y})\big]^{-1} S_N(\mathbf{y}),
		\)
		and hence
		\[
		\sup_{\mathbf{y}\in\mathcal{Y}_0}
		\big\|\hat\gamma(\mathbf{y}) - \gamma^\ast(\mathbf{y})\big\|
		\;\le\;
		\sup_{\mathbf{y}\in\mathcal{Y}_0}
		\big\|\big[H_N^\ast(\mathbf{y})\big]^{-1}\big\|
		\cdot
		\sup_{\mathbf{y}\in\mathcal{Y}_0} \big\|S_N(\mathbf{y})\big\|
		\;\le\;
		\frac{2}{c_H}\,
		\sup_{\mathbf{y}\in\mathcal{Y}_0} \big\|S_N(\mathbf{y})\big\|.
		\]
		
		Using the bound from Lemma~\ref{lem:score-expansion}, we finally obtain
		\[
		\sup_{\mathbf{y}\in\mathcal{Y}_0}
		\big\|\hat\gamma(\mathbf{y}) - \gamma^\ast(\mathbf{y})\big\|
		=
		O_p\!\left( h^{p+1} + \sqrt{\frac{\log (1/h)}{N h^s}} \right),
		\]
		as claimed.
	\end{proof}

	\begin{theorem}\label{thm:uniform-rate-theta}
		Assume \textnormal{(R1)--(R7)} hold with $p$ odd. Let \(\hat\nu(\mathbf{y}) = e_0^\top \hat\gamma(\mathbf{y}), \
		\nu(\mathbf{y}) = e_0^\top \gamma^\ast(\mathbf{y}),
		\)
		and define $\theta(\mathbf{y}) = \psi^{-1}(\nu(\mathbf{y}))$ and
		$\hat\theta(\mathbf{y}) = \psi^{-1}(\hat\nu(\mathbf{y}))$. Then, as $N\to\infty$,
		\begin{align}
			\sup_{\mathbf{y}\in\mathcal{Y}_0} \big|\hat\nu(\mathbf{y}) - \nu(\mathbf{y})\big|
			&=
			O_p\!\left( h^{p+1} + \sqrt{\frac{\log (1/h)}{N h^s}} \right),
			\label{eq:nu-uniform-rate}
			\\[0.25cm]
			\sup_{\mathbf{y}\in\mathcal{Y}_0} \big|\hat\theta(\mathbf{y}) - \theta(\mathbf{y})\big|
			&=
			O_p\!\left( h^{p+1} + \sqrt{\frac{\log (1/h)}{N h^s}} \right).
			\label{eq:theta-uniform-rate}
		\end{align}
		In particular, under \textnormal{(R3)}, both $\hat\nu(\cdot)$ and
		$\hat\theta(\cdot)$ are uniformly consistent on $\mathcal{Y}_0$. Applying Remark~\ref{rem:logN-rate} we can replace $\log(1/h)$ with $\log N$ for sufficiently large $N$.
	\end{theorem}
	
	\begin{proof}
		The result for $\hat\nu(\cdot)$ is an immediate consequence of the uniform rate for
		$\hat\gamma(\cdot)$ in Lemma~\ref{lem:gamma-uniform-rate}. Recall that
		$\nu(\mathbf{y}) = e_0^\top \gamma^\ast(\mathbf{y})$ and $\hat\nu(\mathbf{y}) = e_0^\top \hat\gamma(\mathbf{y})$, where
		$e_0 = (1,0,\dots,0)^\top \in \mathbb{R}^{d_p}$, so that
		\[
		\hat\nu(\mathbf{y}) - \nu(\mathbf{y})
		=
		e_0^\top \bigl(\hat\gamma(\mathbf{y}) - \gamma^\ast(\mathbf{y})\bigr).
		\]
		Hence, for each $\mathbf{y} \in \mathcal{Y}_0$,
		\[
		\bigl|\hat\nu(\mathbf{y}) - \nu(\mathbf{y})\bigr|
		\le
		\|e_0\|\,
		\bigl\|\hat\gamma(\mathbf{y}) - \gamma^\ast(\mathbf{y})\bigr\|
		=
		\bigl\|\hat\gamma(\mathbf{y}) - \gamma^\ast(\mathbf{y})\bigr\|,
		\]
		and taking suprema in $\mathbf{y}$ gives
		\[
		\sup_{\mathbf{y}\in\mathcal{Y}_0}
		\bigl|\hat\nu(\mathbf{y}) - \nu(\mathbf{y})\bigr|
		\le
		\sup_{\mathbf{y}\in\mathcal{Y}_0}
		\bigl\|\hat\gamma(\mathbf{y}) - \gamma^\ast(\mathbf{y})\bigr\|.
		\]
		Thus, by Lemma~\ref{lem:gamma-uniform-rate},
		\[
		\sup_{\mathbf{y}\in\mathcal{Y}_0}
		\bigl|\hat\nu(\mathbf{y}) - \nu(\mathbf{y})\bigr|
		=
		O_p\!\left(
		h^{p+1}
		+
		\sqrt{\frac{\log(1/h)}{N h^s}}
		\right),
		\]
		which proves \eqref{eq:nu-uniform-rate}.
		
		\medskip\noindent
		We now turn to $\hat\theta(\cdot)$. By definition,
		\[
		\theta(\mathbf{y}) = \psi^{-1}\bigl(\nu(\mathbf{y})\bigr),
		\qquad
		\hat\theta(\mathbf{y}) = \psi^{-1}\bigl(\hat\nu(\mathbf{y})\bigr).
		\]
		For each fixed $y \in \mathcal{Y}_0$, the mean value theorem yields
		\[
		\hat\theta(\mathbf{y}) - \theta(\mathbf{y})
		=
		\bigl(\psi^{-1}\bigr)'\!\bigl(\xi_N(\mathbf{y})\bigr)\,
		\bigl(\hat\nu(\mathbf{y}) - \nu(\mathbf{y})\bigr),
		\]
		for some intermediate point $\xi_N(\mathbf{y})$ lying between $\hat\nu(\mathbf{y})$ and $\nu(\mathbf{y})$.
		
		By (R6), $\psi^{-1}$ is continuously differentiable and its derivative
		$\bigl(\psi^{-1}\bigr)'$ is uniformly bounded on the compact set
		$\{\nu(\mathbf{y}): \mathbf{y}\in\mathcal{Y}_0\}$. Moreover, by \eqref{eq:nu-uniform-rate} and the fact that
		$\mathcal{Y}_0$ is compact, we have
		\[
		\sup_{\mathbf{y}\in\mathcal{Y}_0}
		\bigl|\hat\nu(\mathbf{y}) - \nu(\mathbf{y})\bigr| \xrightarrow{p} 0,
		\]
		so for all sufficiently large $N$, the random values $\hat\nu(\mathbf{y})$ lie, with probability
		tending to one, in a fixed compact neighborhood of $\{\nu(\mathbf{y}):\mathbf{y}\in\mathcal{Y}_0\}$. By
		continuity of $\bigl(\psi^{-1}\bigr)'$, there exists a constant $C_\psi < \infty$ such that
		\[
		\sup_{\mathbf{y}\in\mathcal{Y}_0}
		\bigl| \bigl(\psi^{-1}\bigr)'(\xi_N(\mathbf{y})) \bigr|
		\le
		C_\psi
		\quad\text{with probability tending to one.}
		\]
		
		Consequently,
		\[
		\sup_{\mathbf{y}\in\mathcal{Y}_0}
		\bigl|\hat\theta(\mathbf{y}) - \theta(\mathbf{y})\bigr|
		\le
		C_\psi
		\sup_{\mathbf{y}\in\mathcal{Y}_0}
		\bigl|\hat\nu(\mathbf{y}) - \nu(\mathbf{y})\bigr|
		=
		O_p\!\left(
		h^{p+1}
		+
		\sqrt{\frac{\log(1/h)}{N h^s}}
		\right),
		\]
		which establishes \eqref{eq:theta-uniform-rate}.
		
		\medskip\noindent
		Finally, under (R3) we have $h_N \to 0$ and
		$N h_N^s / \log(1/h_N) \to \infty$, so both $h^{p+1}$ and
		$\sqrt{\log(1/h)/(N h^s)}$ converge to zero. Thus the right-hand sides of
		\eqref{eq:nu-uniform-rate} and \eqref{eq:theta-uniform-rate} converge to zero, and both
		$\hat\nu(\cdot)$ and $\hat\theta(\cdot)$ are uniformly consistent on $\mathcal{Y}_0$.
	\end{proof}

	\begin{corollary}\label{cor:tau-uniform-rate}
		Let $\tau(\theta)$ denote the Kendall's $\tau$ associated with the chosen copula
		family, and define $\tau(\mathbf{y}) = \tau(\theta(\mathbf{y}))$ and $\hat\tau(\mathbf{y}) = \tau(\hat\theta(\mathbf{y}))$.
		Under the assumptions of Theorem~\ref{thm:uniform-rate-theta}, suppose in addition
		that $\tau(\cdot)$ is continuously differentiable on $\Theta_0$.
		Then, as $N\to\infty$,
		\[
		\sup_{\mathbf{y}\in\mathcal{Y}_0} \big|\hat\tau(\mathbf{y}) - \tau(\mathbf{y})\big|
		=
		O_p\!\left( h^{p+1} + \sqrt{\frac{\log (1/h)}{N h^s}} \right),
		\]
		and in particular $\hat\tau(\cdot)$ is uniformly consistent on $\mathcal{Y}_0$. Applying Remark~\ref{rem:logN-rate} we can replace $\log(1/h)$ with $\log N$ for sufficiently large $N$.
	\end{corollary}
	
	\begin{proof}
		By Theorem~\ref{thm:uniform-rate-theta}, under \textnormal{(R1)--(R7)} with $p$ odd we have
		\[
		\sup_{\mathbf{y}\in\mathcal{Y}_0} \big|\hat\theta(\mathbf{y}) - \theta(\mathbf{y})\big|
		=
		O_p\!\left( h^{p+1} + \sqrt{\frac{\log (1/h)}{N h^s}} \right).
		\]
		
		For each fixed $\mathbf{y}\in\mathcal{Y}_0$, the mean value theorem applied to the
		one-dimensional map $\tau(\cdot)$ gives
		\[
		\hat\tau(\mathbf{y}) - \tau(\mathbf{y})
		=
		\tau\big(\hat\theta(\mathbf{y})\big) - \tau\big(\theta(\mathbf{y})\big)
		=
		\tau'\big(\xi_N(\mathbf{y})\big)\,\big(\hat\theta(\mathbf{y})-\theta(\mathbf{y})\big),
		\]
		for some intermediate point $\xi_N(\mathbf{y})$ between $\hat\theta(\mathbf{y})$ and
		$\theta(\mathbf{y})$.
		
		By assumption, $\tau(\cdot)$ is continuously differentiable on $\Theta_0$.
		Moreover, by Theorem~\ref{thm:uniform-rate-theta} and the compactness of
		$\mathcal{Y}_0$, we have
		\[
		\sup_{\mathbf{y}\in\mathcal{Y}_0} \big|\hat\theta(\mathbf{y}) - \theta(\mathbf{y})\big| \xrightarrow{p} 0,
		\]
		so for all sufficiently large $N$ the random values $\hat\theta(\mathbf{y})$ lie, with
		probability tending to one, in any fixed compact neighborhood of
		$\{\theta(\mathbf{y}):\mathbf{y}\in\mathcal{Y}_0\}\subset \Theta_0$. By continuity of $\tau'$ on
		$\Theta_0$, there exists a finite constant $C_\tau$ such that
		\[
		\sup_{\mathbf{y}\in\mathcal{Y}_0}
		\big| \tau'\big(\xi_N(\mathbf{y})\big) \big|
		\;\le\;
		C_\tau
		\quad\text{with probability tending to one.}
		\]
		
		Therefore,
		\[
		\sup_{\mathbf{y}\in\mathcal{Y}_0} \big|\hat\tau(\mathbf{y}) - \tau(\mathbf{y})\big|
		\;\le\;
		\sup_{\mathbf{y}\in\mathcal{Y}_0}
		\big| \tau'\big(\xi_N(\mathbf{y})\big) \big|\,
		\sup_{\mathbf{y}\in\mathcal{Y}_0} \big|\hat\theta(\mathbf{y}) - \theta(\mathbf{y})\big|
		\;\le\;
		C_\tau
		\sup_{\mathbf{y}\in\mathcal{Y}_0} \big|\hat\theta(\mathbf{y}) - \theta(\mathbf{y})\big|,
		\]
		and hence
		\[
		\sup_{\mathbf{y}\in\mathcal{Y}_0} \big|\hat\tau(\mathbf{y}) - \tau(\mathbf{y})\big|
		=
		O_p\!\left( h^{p+1} + \sqrt{\frac{\log (1/h)}{N h^s}} \right).
		\]
		
		Under the bandwidth condition \textnormal{(R3)}, both $h^{p+1}$ and
		$\sqrt{\log(1/h)/(N h^s)}$ converge to zero, so the right-hand side converges
		to zero, and $\hat\tau(\cdot)$ is uniformly consistent on $\mathcal{Y}_0$.
	\end{proof}
	
\section{Numerical Analysis}\label{sec:simulation}

We study the finite-sample performance of the proposed local likelihood estimator and assess the empirical relevance of the uniform convergence rate derived in Section~\ref{sec:main-results}. Particular emphasis is placed on evaluating whether the stochastic term involving $\log(1/h)$ can be accurately approximated by $\log N$ in practice.

\paragraph{Data-generating mechanism.}
The covariate $\mathbf{Y}$ is generated from a truncated normal distribution
$\mathrm{TN}(0,4;[-2,2])$. Conditional on $Y=y$ $(s=1)$, the pair $(U_1,U_2)$ follows a Frank copula with parameter $\theta(y)=\psi^{-1}(\nu(y))$, where the inverse link function is
\(
\psi^{-1}(\nu)=1+\frac{4}{1+e^{-\nu}},
\)
so that $\theta(y)\in(1,5)$. Estimation is carried out using local linear likelihood smoothing with $p=1$ and the Epanechnikov kernel. Maximization of the local log-likelihood is performed using a BFGS algorithm. The bandwidth $h$ is selected by least-squares cross-validation via maximization of~\eqref{lcv_h}, without imposing an upper bound.

\paragraph{Models.}
Two calibration functions are considered:
\begin{align}
	\theta_{\mathrm{M1}}(y) &= 5-y^2, \label{eq:thetaM1}\\
	\theta_{\mathrm{M2}}(y) &= 3+y+\frac{5}{6}\Big(e^{-y^2}-e^{-4}\Big). \label{eq:thetaM2}
\end{align}
Model~M1 exhibits strong curvature and substantial spatial variation, while Model~M2 represents a smoother nonlinear trend with more moderate local variation. These two designs allow us to examine the robustness of the uniform rate approximation across different functional shapes.

\begin{figure}[ht]
	\centering
	\begin{subfigure}[t]{0.48\textwidth}
		\centering
		\includegraphics[width=\textwidth]{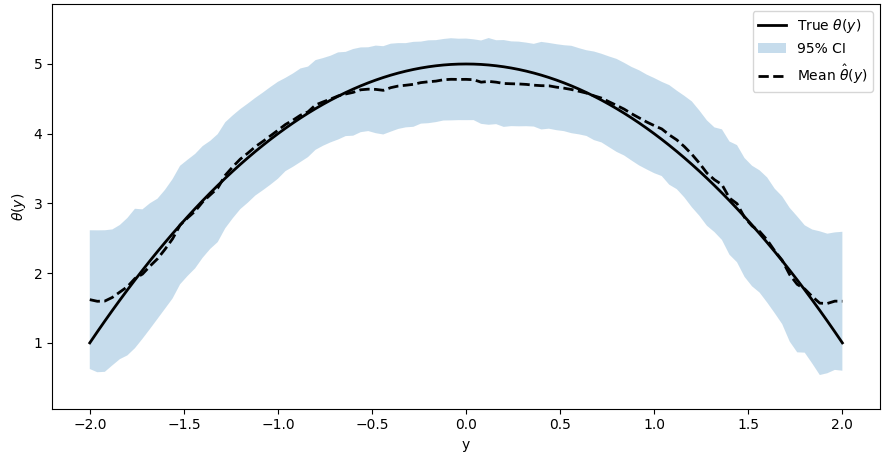}
		\caption{Model M1 (Frank copula).}
	\end{subfigure}
	\hfill
	\begin{subfigure}[t]{0.48\textwidth}
		\centering
		\includegraphics[width=\textwidth]{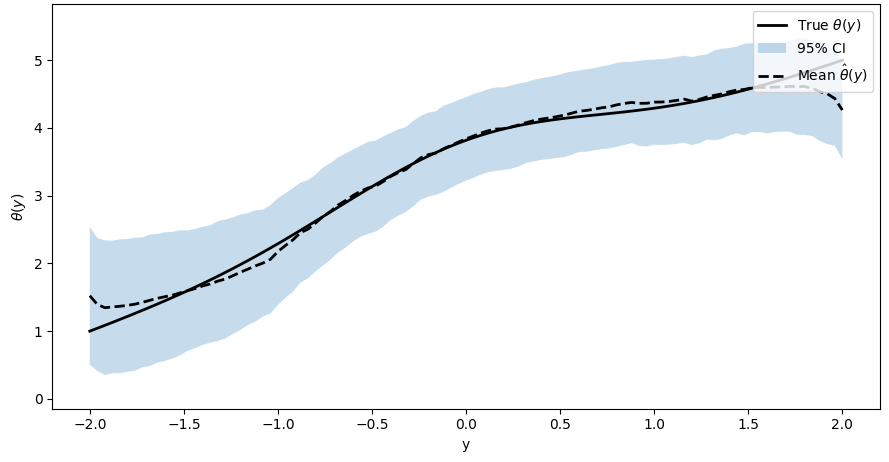}
		\caption{Model M2 (Frank copula).}
	\end{subfigure}
	\caption{\it True parameter curve $\theta(y)$ (solid), Monte Carlo mean of $\widehat{\theta}(y)$ (dashed), and pointwise $95\%$ Monte Carlo bands across replications (shaded), shown for a representative sample size.}
	\label{fig:thetafit-M1M2}
\end{figure}

\paragraph{Illustration of functional estimation.}
Figure~\ref{fig:thetafit-M1M2} compares the true function $\theta(\cdot)$ with the Monte Carlo mean of the estimator $\widehat{\theta}(\cdot)$ for a representative sample size. Pointwise $95\%$ Monte Carlo bands are also shown. The results are based on samples of size $N=500$ and $100$ Monte Carlo replications. In both models, the estimator closely follows the global shape of $\theta(\cdot)$. As expected, variability increases near the boundaries of the covariate domain, where effective local sample sizes are smaller and kernel smoothing becomes less stable. These boundary regions are excluded when computing the uniform error reported below. Further details on the bias, variance, and asymptotic normality of $\widehat{\theta}(\cdot)$ and $\widehat{\boldsymbol{\gamma}}(\cdot)$ in~\eqref{estimator_main} are provided in the Appendix.

\paragraph{Performance metric and rate proxy.}
For each Monte Carlo replication, performance is summarized by the uniform estimation error
\(
\sup_{\mathbf{y}\in\mathcal{Y}_0}\bigl|\widehat{\theta}(y)-\theta(y)\bigr|,
\quad
\mathcal{Y}_0=[-2+\delta,\,2-\delta], \ \delta=0.2,
\)
which reduces boundary effects. To compare empirical performance with the theoretical results, we also compute the rate proxy
\(
h^{p+1}+\sqrt{\frac{\log N}{N h^{s}}},
\quad p=s=1,
\)
using the same bandwidth $h$ selected by cross-validation. As discussed in Section~\ref{sec:main-results}, this quantity captures the dominant uniform stochastic fluctuations when $\log N$ replaces $\log(1/h)$; see Remark~\ref{rem:logN-rate}.

Simulations are conducted for sample sizes $N\in\{100,200,500,750,1000,1500,2000\}$, with $R=100$ replications per configuration. For compactness and direct comparison, results for both models are reported together.

\paragraph{Main findings.}
Table~\ref{tab:main-sup-rate} reports the empirical mean, standard deviation, and median of the uniform estimation error (\textit{sup}) and of the corresponding rate proxy (\textit{rate}). For both models, the uniform error decreases as the sample size increases, and the rate proxy exhibits a closely aligned decline. Agreement between the empirical error and the theoretical rate is particularly strong for moderate and large sample sizes, supporting the practical replacement of $\log(1/h)$ by $\log N$, even when cross-validation occasionally selects relatively large bandwidths for smaller samples.

\begin{table}[!ht]
	\centering
	\small
	\setlength{\tabcolsep}{6pt}
	\renewcommand{\arraystretch}{1.05}
	
	\begin{tabular*}{\textwidth}{@{\extracolsep{\fill}} l r r r r r r r}
		\toprule
		Model & $N$ &
		$\mathrm{sup}_{\mathrm{mean}}$ &
		$\mathrm{sup}_{\mathrm{sd}}$ &
		$\mathrm{sup}_{\mathrm{median}}$ &
		$\mathrm{rate}_{\mathrm{mean}}$ &
		$\mathrm{rate}_{\mathrm{sd}}$ &
		$\mathrm{rate}_{\mathrm{median}}$ \\
		\midrule
		\multirow{7}{*}{M1}
		& 100  & 2.4138 & 0.7896 & 2.5546 & 1.4981 & 0.8618 & 1.3600 \\
		& 200  & 1.8559 & 0.8405 & 1.9266 & 1.3065 & 0.6611 & 1.5206 \\
		& 500  & 0.9964 & 0.4890 & 0.9022 & 0.9802 & 0.4305 & 1.0788 \\
		& 750  & 0.9342 & 0.5532 & 0.7647 & 0.8096 & 0.3839 & 0.8022 \\
		& 1000 & 0.7907 & 0.4390 & 0.6666 & 0.7421 & 0.3385 & 0.8291 \\
		& 1500 & 0.6189 & 0.3373 & 0.5364 & 0.6483 & 0.2862 & 0.7090 \\
		& 2000 & 0.5752 & 0.2485 & 0.5313 & 0.5990 & 0.2332 & 0.6337 \\
		\midrule
		\multirow{7}{*}{M2}
		& 100  & 2.0354 & 1.0721 & 1.9338 & 1.7376 & 0.8866 & 2.0050 \\
		& 200  & 1.3672 & 0.8389 & 1.1071 & 1.3800 & 0.6337 & 1.5371 \\
		& 500  & 0.8417 & 0.5490 & 0.7144 & 1.1297 & 0.4066 & 1.4144 \\
		& 750  & 0.7060 & 0.5676 & 0.5560 & 0.9808 & 0.3281 & 1.2102 \\
		& 1000 & 0.6776 & 0.5175 & 0.5370 & 0.8176 & 0.3185 & 1.0783 \\
		& 1500 & 0.5601 & 0.3242 & 0.4595 & 0.7372 & 0.2664 & 0.9269 \\
		& 2000 & 0.4581 & 0.1902 & 0.4285 & 0.6803 & 0.2091 & 0.8277 \\
		\bottomrule
	\end{tabular*}
	
	\caption{\it Uniform estimation error
		$\sup_{\mathbf{y}\in\mathcal{Y}_0}\lvert\widehat\theta(\mathbf y)-\theta(\mathbf y)\rvert$
		and theoretical rate proxy
		$h^2+\sqrt{\log N/(Nh)}$ across $R=100$ replications.}
	\label{tab:main-sup-rate}
\end{table}

To quantify scale agreement, Table~\ref{tab:ratio-sup-rate} reports ratios of Monte Carlo means. Ratios close to one indicate that the proxy captures the correct order and scale. For both models, the ratios stabilize with $N$, consistent with the theoretical proxy describing the leading uniform fluctuations. In particular, the ratio behavior is compatible with the theory even though bandwidths may occasionally exceed one in smaller samples (see below).

\begin{table}[!ht]
	\centering
	\small
	\setlength{\tabcolsep}{8pt}
	\renewcommand{\arraystretch}{1.05}
	
	\begin{tabular*}{\textwidth}{@{\extracolsep{\fill}} l r r r}
		\toprule
		Model & $N$ &
		$\mathrm{sup}_{\mathrm{mean}}/\mathrm{rate}_{\mathrm{mean}}$ &
		$\mathrm{rate}_{\mathrm{mean}}/\mathrm{sup}_{\mathrm{mean}}$ \\
		\midrule
		\multirow{7}{*}{M1}
		& 100  & 1.6112 & 0.6207 \\
		& 200  & 1.4205 & 0.7040 \\
		& 500  & 1.0165 & 0.9837 \\
		& 750  & 1.1539 & 0.8666 \\
		& 1000 & 1.0655 & 0.9386 \\
		& 1500 & 0.9546 & 1.0475 \\
		& 2000 & 0.9602 & 1.0414 \\
		\midrule
		\multirow{7}{*}{M2}
		& 100  & 1.1724 & 0.8533 \\
		& 200  & 0.9901 & 1.0100 \\
		& 500  & 0.7456 & 1.3410 \\
		& 750  & 0.7204 & 1.3897 \\
		& 1000 & 0.8290 & 1.2067 \\
		& 1500 & 0.7606 & 1.3162 \\
		& 2000 & 0.6748 & 1.4851 \\
		\bottomrule
	\end{tabular*}
	
	\caption{\it Ratios comparing Monte Carlo means of the uniform estimation error and the theoretical rate.}
	\label{tab:ratio-sup-rate}
\end{table}

\paragraph{Bandwidth behavior (cross-validation).}
Table~\ref{tab:bandwidth-summaries} summarizes the distribution of selected bandwidths. For small and moderate $N$, cross-validation can prefer global smoothing and occasionally yields $h>1$. As $N$ increases, the distribution shifts downward and concentrates, consistent with $h\to0$ in probability. Importantly, the presence of some $h>1$ values does not contradict the theoretical stochastic term used in the paper: under the entropy control in Lemma~\ref{lem:kernel-class-entropy} and Remark~\ref{rem:logN-rate}, the uniform bound is expressed with $\log N$ rather than $\log(1/h)$, making the simulation design directly relevant to the theory being assessed.

\begin{table}[!ht]
	\centering
	\small
	\setlength{\tabcolsep}{8pt}
	\renewcommand{\arraystretch}{1.05}
	
	\begin{tabular*}{\textwidth}{@{\extracolsep{\fill}} l r r r r r r}
		\toprule
		Model & $N$ &
		$h_{\mathrm{mean}}$ &
		$h_{\mathrm{sd}}$ &
		$h_{\mathrm{median}}$ &
		$h_{\mathrm{min}}$ &
		$h_{\mathrm{max}}$ \\
		\midrule
		\multirow{7}{*}{M1}
		& 100  & 1.0145 & 0.4780 & 1.0737 & 0.2272 & 1.5875 \\
		& 200  & 0.9791 & 0.4069 & 1.1705 & 0.2029 & 1.3856 \\
		& 500  & 0.8813 & 0.2809 & 0.9830 & 0.1709 & 1.1534 \\
		& 750  & 0.7824 & 0.2883 & 0.8331 & 0.1576 & 1.0637 \\
		& 1000 & 0.7580 & 0.2586 & 0.8599 & 0.2904 & 1.0047 \\
		& 1500 & 0.7100 & 0.2381 & 0.7942 & 0.1388 & 0.9264 \\
		& 2000 & 0.6939 & 0.2007 & 0.7500 & 0.1312 & 0.8744 \\
		\midrule
		\multirow{7}{*}{M2}
		& 100  & 1.1382 & 0.4694 & 1.3499 & 0.2181 & 1.5871 \\
		& 200  & 1.0312 & 0.3793 & 1.1785 & 0.2030 & 1.3863 \\
		& 500  & 0.9644 & 0.2802 & 1.1458 & 0.1719 & 1.1530 \\
		& 750  & 0.9043 & 0.2431 & 1.0588 & 0.1589 & 1.0644 \\
		& 1000 & 0.8090 & 0.2529 & 0.9984 & 0.1507 & 1.0053 \\
		& 1500 & 0.7702 & 0.2400 & 0.9243 & 0.1396 & 0.9267 \\
		& 2000 & 0.7601 & 0.1747 & 0.8731 & 0.1311 & 0.8752 \\
		\bottomrule
	\end{tabular*}
	
	\caption{\it Summaries of the cross-validated bandwidth $h$ across $R=100$ replications.}
	\label{tab:bandwidth-summaries}
\end{table}

\begin{figure}[ht]
		\centering
		\begin{subfigure}[t]{0.48\textwidth}
			\centering
			\includegraphics[width=\textwidth]{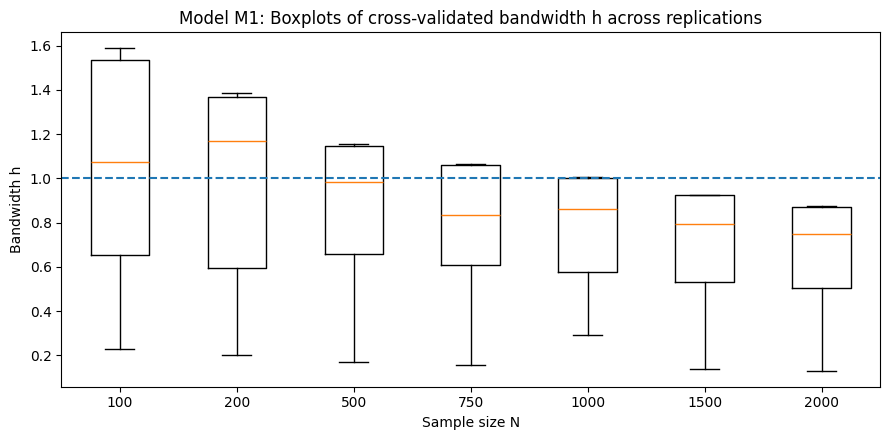}
			\caption{Model M1.}
		\end{subfigure}
		\hfill
		\begin{subfigure}[t]{0.48\textwidth}
			\centering
			\includegraphics[width=\textwidth]{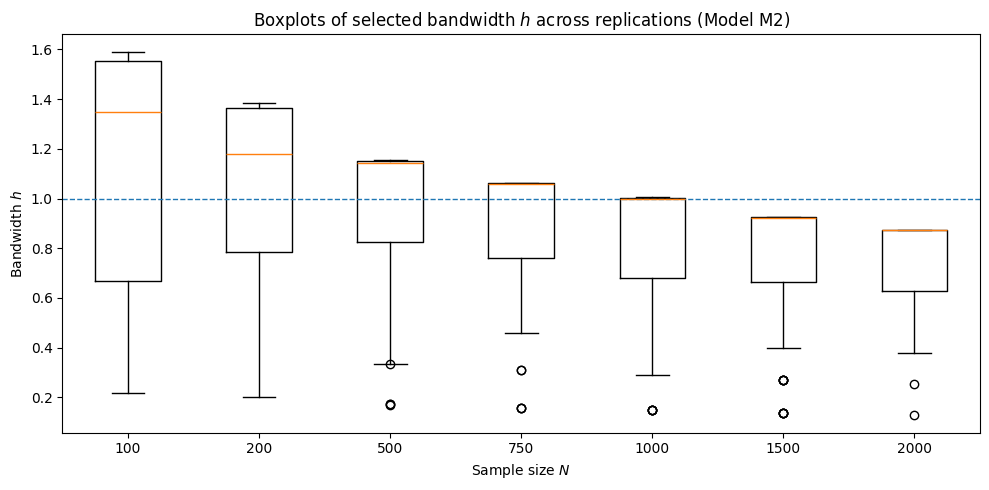}
			\caption{Model M2.}
		\end{subfigure}
		\caption{\it Boxplots of the cross-validated bandwidth $h$ across $R=100$ replications for each sample size. The dashed horizontal line marks $h=1$.}
		\label{fig:bandwidth-M1M2}
	\end{figure}

	Figure~\ref{fig:bandwidth-M1M2} provides the corresponding boxplots side-by-side for Models M1 and M2. The downward shift and tightening of the bandwidth distribution with $N$ is clearly visible in both designs.
	
\begin{figure}[ht]
		\centering
		\begin{subfigure}[t]{0.48\textwidth}
			\centering
			\includegraphics[width=\textwidth]{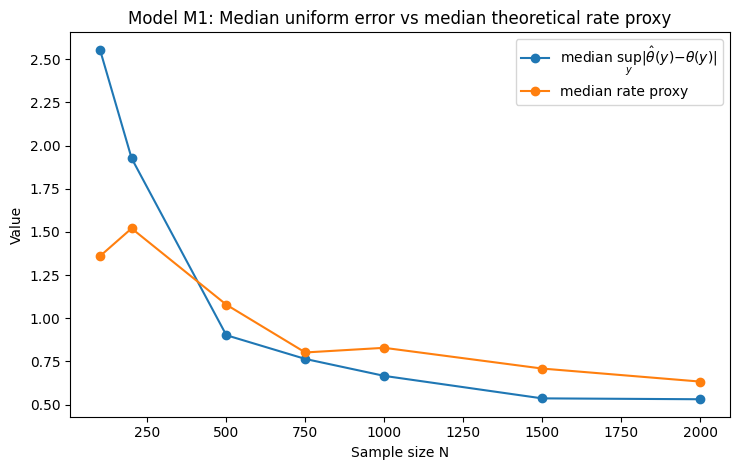}
			\caption{Model M1.}
		\end{subfigure}
		\hfill
		\begin{subfigure}[t]{0.49\textwidth}
			\centering
			\includegraphics[width=\textwidth]{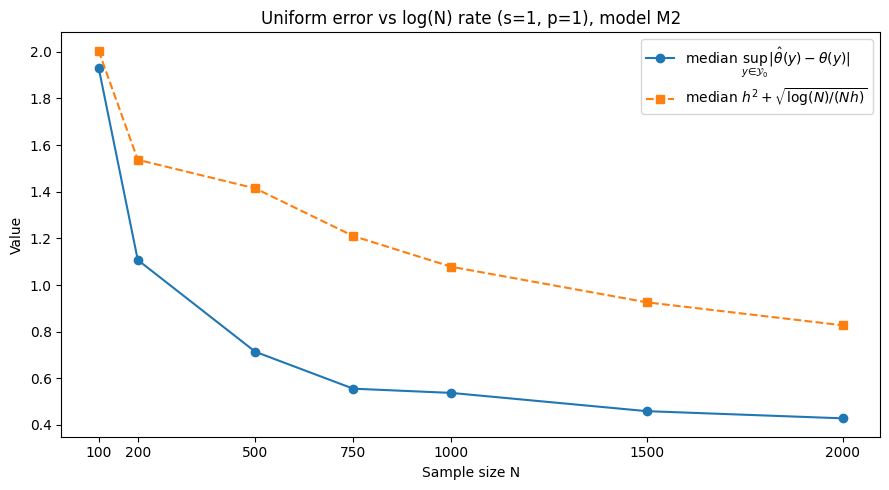}
			\caption{Model M2.}
		\end{subfigure}
		\caption{\it Median uniform estimation error and median theoretical rate $\,h^2+\sqrt{\log N/(Nh)}\,$ as functions of $N$.}
		\label{fig:ratecurves-M1M2}
	\end{figure}

	\paragraph{Median trends (uniform error vs.\ rate).}
	Figure~\ref{fig:ratecurves-M1M2} compares the median uniform error to the median rate as functions of $N$, shown side-by-side for the two models. The two curves move in tandem in both cases, which is the key empirical check: the rate proxy captures not only the order but also the practical decay pattern of the uniform error when $\log N$ is used.
	
\begin{figure}[ht]
		\centering
		\begin{subfigure}[t]{0.99\textwidth}
			\centering
			\includegraphics[width=\textwidth]{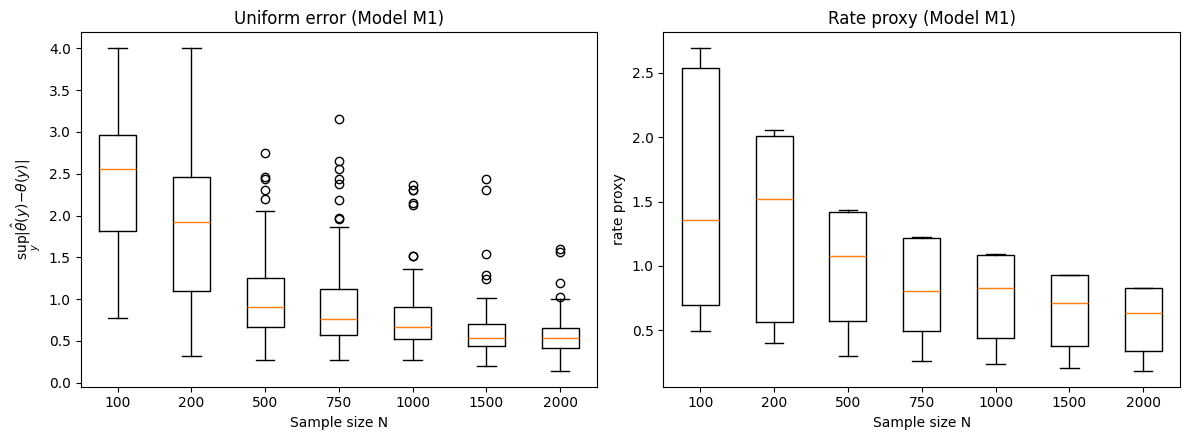}
			\caption{Model M1: uniform error and rate.}
		\end{subfigure}

		\vspace{0.35cm}
		
		\begin{subfigure}[t]{0.48\textwidth}
			\centering
			\includegraphics[width=\textwidth]{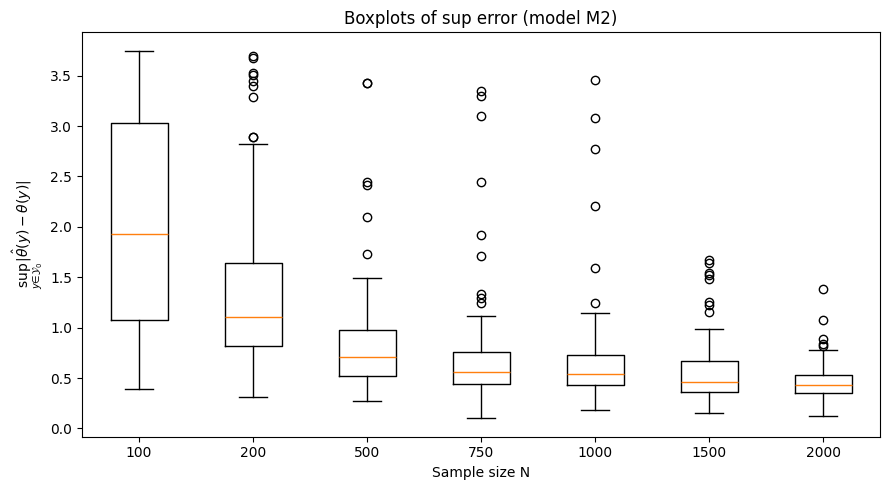}
			\caption{Model M2: uniform error.}
		\end{subfigure}
		\hfill
		\begin{subfigure}[t]{0.48\textwidth}
			\centering
			\includegraphics[width=\textwidth]{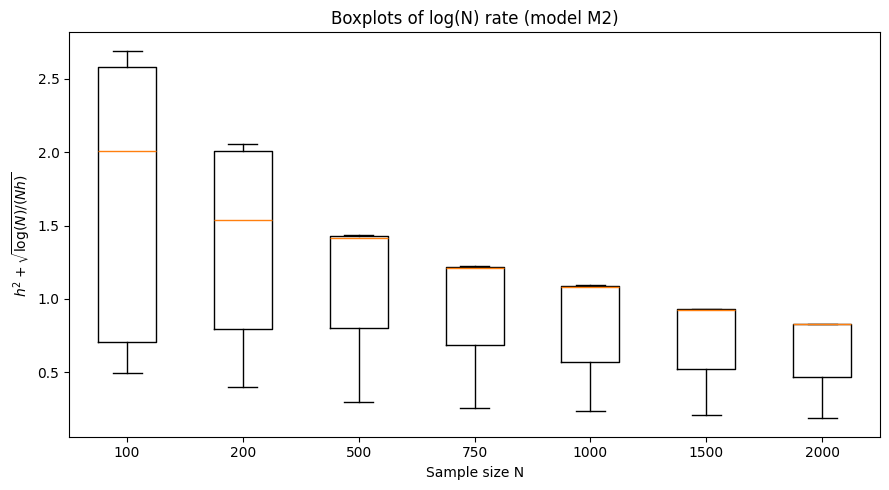}
			\caption{Model M2: rate.}
		\end{subfigure}
		\caption{\it Boxplots of the uniform estimation error and the theoretical rate $\,h^2+\sqrt{\log N/(Nh)}\,$ across $R=100$ replications for each sample size, shown as paired panels for Model M1 (top row) and Model M2 (bottom row).}
		\label{fig:boxpairs-M1M2}
	\end{figure}

	\paragraph{Distributional comparison across replications.}
	Finally, Figure~\ref{fig:boxpairs-M1M2} reports boxplots of the uniform estimation error and the theoretical rate across replications, arranged as paired panels (error vs.\ rate) for Model M1 and then Model M2. In both models, dispersion contracts as $N$ grows. The rate distribution mirrors the error distribution in its contraction and its scale evolution, providing a stringent confirmation of the uniform-rate approximation under realistic, data-driven bandwidths.
	
\section{Discussion}

This paper develops uniform consistency results for local likelihood estimators of covariate-dependent copula parameters in multivariate settings. In contrast to much of the existing literature, which concentrates on pointwise asymptotic behavior at fixed covariate values, the analysis here provides uniform control over compact regions of the covariate space. This distinction is substantive rather than purely technical. Uniform convergence is required to ensure the reliability of local optimization procedures, to justify inference that holds simultaneously across covariate values, and to support bandwidth selection strategies that depend on global rather than pointwise performance.

From a methodological standpoint, the main contribution is the derivation of explicit uniform convergence rates for the local log-likelihood and its first two derivatives. These results imply that the local maximum likelihood estimator converges uniformly to its population target, and that the estimated calibration and copula parameter functions inherit the same rates. Uniform control of the Hessian plays a central role, as it guarantees the existence and stability of local maximizers with high probability and ensures that Newton–Raphson or quasi-Newton algorithms behave well throughout the covariate domain. This provides theoretical support for the widespread empirical use of local likelihood methods in conditional copula modeling.

The theoretical arguments rely on empirical process techniques for kernel-weighted function classes indexed jointly by covariate locations and local polynomial coefficients. A principal difficulty arises from kernel localization, which produces shrinking neighborhoods and envelopes that increase as the bandwidth decreases. By combining polynomial entropy bounds with geometric properties of the covariate space, these challenges can be addressed under relatively mild smoothness and moment assumptions. The resulting convergence rates explicitly depend on the dimension of the covariate, reflecting the familiar curse of dimensionality in multivariate smoothing.

The results obtained here extend earlier work on local likelihood estimation for conditional copula models. While pointwise bias and variance expansions remain informative for local behavior, uniform convergence is essential for establishing global consistency of dependence measures such as Kendall’s $\tau$ and for validating inference procedures based on supremum-type error bounds. The uniform consistency theory developed in this paper therefore provides a natural foundation for constructing uniform confidence bands, conducting simultaneous hypothesis tests, and developing other global inferential tools in covariate-dependent copula models.

Several directions for future research remain open. Although the present results offer theoretical justification for data-driven bandwidth selection, a complete asymptotic analysis of cross-validated bandwidths in conditional copula models has yet to be developed. In addition, the current framework assumes independent observations. Extending the uniform results to dependent settings, including time series and spatial data, would substantially broaden their applicability. Finally, while the focus here is on single-parameter copula families, many applications involve multi-parameter copulas or vine constructions, where uniform control may require additional structure.

The empirical process techniques used in this paper are not restricted to copula models and may be useful in a wider class of kernel-weighted $M$-estimation problems with covariate-indexed parameters. As a result, the methods and results developed here may find applications beyond conditional copula models, particularly in semiparametric frameworks that combine local smoothing with global constraints.

	\appendix
	\section{Appendix}
	\subsection{Auxiliary multivariate local likelihood results}
	
	This appendix collects pointwise asymptotic results for multivariate local likelihood
	estimators that complement the uniform convergence theory developed in the main text.
	The results summarized here extend standard local polynomial likelihood arguments
	(Fan and Gijbels \cite{fan1996local}) to the conditional copula setting with link-transformed parameters
	and provide explicit bias, variance, and asymptotic normality statements (see Acar et al. \cite{acar2011dependence}, Muia et al. \cite{muia2025local} and  \cite{muia2025multivariate}).
	They also clarify the origin of the smoothing and stochastic terms that appear in the
	uniform rates established in Section~4.
	
	Throughout, let $\mathbf Y\in\mathbb R^s$ be a covariate vector with density
	$f_{\mathbf Y}$, and let $\nu(\mathbf y)$ be a real-valued calibration function.
	The copula parameter is given by
	\(\theta(\mathbf y)=\psi^{-1}\!\{\nu(\mathbf y)\},\)
	where $\psi^{-1}$ is a strictly monotone inverse link.
	
	Let $\alpha\in\mathbb N^s$ be a multi-index, with
	\(
	|\alpha| = \sum_{j=1}^s \alpha_j, \quad
	\alpha! = \prod_{j=1}^s \alpha_j!, \quad
	\mathbf v^\alpha = \prod_{j=1}^s v_j^{\alpha_j},
	D^\alpha
	=
	\frac{\partial^{|\alpha|}}
	{\partial y_1^{\alpha_1}\cdots\partial y_s^{\alpha_s}}.
	\)
	We approximate $\nu(\mathbf Y_i)$ locally by
	\(\nu(\mathbf Y_i)\approx
	\boldsymbol\phi_p(\mathbf Y_i-\mathbf y)^\top\boldsymbol\gamma,\)
	and estimate the derivative by
	\(
	\widehat{D^\alpha\nu}(\mathbf y)
	=
	\alpha!\,\hat\gamma_\alpha,
	\quad |\alpha|\le p,
	\)
	where $\hat\gamma_\alpha$ denotes the coefficient associated with $\alpha$ and
	$\mathbf e_\alpha$ selects this coordinate in $\boldsymbol\gamma$.
	
	We consider a local polynomial likelihood fit of degree $p$ using an isotropic
	bandwidth matrix $H=h\mathbf I_s$.
	Since $H$ is diagonal with identical entries, smoothing is applied uniformly
	along the coordinate axes of the covariates.
	
	Define the kernel moment matrices
	\(
	S=
	\int
	\boldsymbol\phi_p(\mathbf v)\boldsymbol\phi_p(\mathbf v)^\top
	K(\mathbf v)\,d\mathbf v,
	\
	S^*
	=
	\int
	\boldsymbol\phi_p(\mathbf v)\boldsymbol\phi_p(\mathbf v)^\top
	K^2(\mathbf v)\,d\mathbf v,
	\)
	and the kernel moment vectors
	\(
	\mathbf c_p
	=
	\int
	\boldsymbol\phi_p(\mathbf v)\mathbf r_{p+1}(\mathbf v)
	K(\mathbf v)\,d\mathbf v,
	\quad
	\tilde{\mathbf c}_p
	=
	\int
	\boldsymbol\phi_p(\mathbf v)\mathbf r_{p+2}(\mathbf v)
	K(\mathbf v)\,d\mathbf v,
	\)
	where $\mathbf r_{p+1}$ and $\mathbf r_{p+2}$ collect monomials of total degree
	$p+1$ and $p+2$, respectively.
	
	Let $\widehat{\boldsymbol\gamma}(\mathbf y)$ denote the local maximizer of the
	kernel-weighted log-likelihood at $\mathbf y$ given in \eqref{estimator_main}.
	Define the conditional Fisher curvature
	\(
	\sigma^2(\mathbf y)
	=
	-\mathbb E\!\left[
	\ell''\!\left(\psi^{-1}(\nu(\mathbf y)),\mathbf U\right)
	\mid \mathbf Y=\mathbf y
	\right],
	\)
	and assume $\sigma^2(\mathbf y)>0$ on a neighborhood of interest.

	\subsection{Asymptotic bias and variance}
	
	\begin{theorem}\label{ThmMulti1}
		Let $f_{\mathbf Y}(\mathbf y)>0$ and assume that $f_{\mathbf Y}$,
		all partial derivatives of $\nu$ up to order $p+2$,
		$\sigma^2(\cdot)$, and $(\psi^{-1})'(\cdot)$
		are continuous in a neighborhood of $\mathbf y$.
		Assume further that $\ell'(\theta,\mathbf U)$ and
		$\ell''(\theta,\mathbf U)$ exist, are continuous in $\theta$,
		and are dominated by integrable envelope functions.
		Let $h\to0$ and $Nh^s\to\infty$ as $N\to\infty$.
		
		\begin{enumerate}
			\item[(i)]
			If $p-|\alpha|$ is odd, then
			\[
			\mathrm{Bias}\{\widehat{D^\alpha\nu}(\mathbf y)\mid\mathbb Y\}
			=
			\alpha!\,h^{\,p+1-|\alpha|}
			\mathbf e_\alpha^\top S^{-1}
			\Bigg(
			\sum_{|\beta|=p+1}
			\frac{D^\beta\nu(\mathbf y)}{\beta!}\,
			\mathbf c_{p,\beta}
			\Bigg)
			+ o_P\!\big(h^{\,p+1-|\alpha|}\big).
			\]
			
			\item[(ii)]
			For any $|\alpha|\le p$,
			\[
			\mathrm{Var}\{\widehat{D^\alpha\nu}(\mathbf y)\mid\mathbb Y\}
			=
			\frac{(\alpha!)^2}{Nh^{\,s+2|\alpha|}
				\,\sigma^2(\mathbf y)\{(\psi^{-1})'(\nu(\mathbf y))\}^2
				f_{\mathbf Y}(\mathbf y)}
			\mathbf e_\alpha^\top S^{-1}S^*S^{-1}\mathbf e_\alpha
			+ o_P\!\Big(\tfrac{1}{Nh^{\,s+2|\alpha|}}\Big).
			\]
		\end{enumerate}
	\end{theorem}
	
	Here $S$ and $S^*$ are kernel moment matrices and
	$\mathbf c_{p,\beta}$ are kernel moment vectors determined by the
	choice of kernel $K$ and basis $\boldsymbol\phi_p$.
	The proof follows standard local likelihood arguments and is provided in Muia et al.~\cite{muia2025multivariate}.
	
	%-------------------------------------------------
	\subsection{Point estimation of the calibration function}
	
	\begin{corollary}\label{cor:nuhat_from_thm}
		Let $p$ be odd and define
		$\widehat\nu(\mathbf y)=\mathbf e_0^\top\widehat{\boldsymbol\gamma}(\mathbf y)$.
		Then
		\begin{align*}
			\mathrm{Bias}\{\widehat\nu(\mathbf y)\mid\mathbb Y\}
			&=
			h^{\,p+1}\,
			\mathbf e_0^\top S^{-1}
			\Bigg(
			\sum_{|\beta|=p+1}
			\frac{D^\beta\nu(\mathbf y)}{\beta!}\,
			\mathbf c_{p,\beta}
			\Bigg)
			+ o_P\!\big(h^{\,p+1}\big),\\[4pt]
			\mathrm{Var}\{\widehat\nu(\mathbf y)\mid\mathbb Y\}
			&=
			\frac{1}{Nh^{\,s}\,
				\sigma^2(\mathbf y)\{(\psi^{-1})'(\nu(\mathbf y))\}^2
				f_{\mathbf Y}(\mathbf y)}
			\mathbf e_0^\top S^{-1}S^*S^{-1}\mathbf e_0
			+ o_P\!\Big(\tfrac{1}{Nh^{\,s}}\Big).
		\end{align*}
	\end{corollary}
	
	%-------------------------------------------------
	\subsection{Copula parameter estimation}
	
	\begin{corollary}\label{cor:thetahat_from_thm}
		Let $\widehat\theta(\mathbf y)=\psi^{-1}\!\{\widehat\nu(\mathbf y)\}$
		and assume $\psi'(\theta(\mathbf y))\neq0$.
		Then
		\begin{align*}
			\mathrm{Bias}\{\widehat\theta(\mathbf y)\mid\mathbb Y\}
			&=
			\frac{1}{\psi'(\theta(\mathbf y))}\;
			h^{\,p+1}\,
			\mathbf e_0^\top S^{-1}
			\Bigg(
			\sum_{|\beta|=p+1}
			\frac{D^\beta\nu(\mathbf y)}{\beta!}\,
			\mathbf c_{p,\beta}
			\Bigg)
			+ o_P\!\big(h^{\,p+1}\big),\\[4pt]
			\mathrm{Var}\{\widehat\theta(\mathbf y)\mid\mathbb Y\}
			&=
			\frac{1}{Nh^{\,s}\,
				\sigma^2(\mathbf y) f_{\mathbf Y}(\mathbf y)}
			\mathbf e_0^\top S^{-1}S^*S^{-1}\mathbf e_0
			+ o_P\!\Big(\tfrac{1}{Nh^{\,s}}\Big).
		\end{align*}
	\end{corollary}
	
	%-------------------------------------------------
	\subsection{Asymptotic normality}
	
	\begin{theorem}[Asymptotic normality]\label{thm:asymp_normality}
		Under the assumptions of Theorem~\ref{ThmMulti1},
		\[
		\sqrt{Nh^{\,s+2|\alpha|}}
		\left(
		\widehat{D^\alpha\nu}(\mathbf y)
		-
		\mathrm{Bias}\{\widehat{D^\alpha\nu}(\mathbf y)\mid\mathbb Y\}
		\right)
		\xrightarrow{d}
		\mathcal N\!\left(
		0,\,
		\frac{(\alpha!)^2 \mathbf e_\alpha^\top S^{-1}S^*S^{-1}\mathbf e_\alpha}{\sigma^2(\mathbf y)\{(\psi^{-1})'(\nu(\mathbf y))\}^2
			f_{\mathbf Y}(\mathbf y)}
		\right).
		\]
	\end{theorem}
	
	\begin{corollary}[Asymptotic normality of $\widehat\nu(\mathbf y)$]
		Let $p$ be odd. Then
		\[
		\sqrt{Nh^{\,s}}
		\left(
		\widehat\nu(\mathbf y)
		-
		\mathrm{Bias}\{\widehat\nu(\mathbf y)\mid\mathbb Y\}
		\right)
		\xrightarrow{d}
		\mathcal N\!\left(
		0,\,
		\frac{\mathbf e_0^\top S^{-1}S^*S^{-1}\mathbf e_0}{\sigma^2(\mathbf y)\{(\psi^{-1})'(\nu(\mathbf y))\}^2
			f_{\mathbf Y}(\mathbf y)}
		\right).
		\]
	\end{corollary}
	
	\begin{corollary}[Asymptotic normality of $\widehat\theta(\mathbf y)$]
		\[
		\sqrt{Nh^{\,s}}
		\left(
		\widehat\theta(\mathbf y)
		-
		\theta(\mathbf y)
		-
		\mathrm{Bias}\{\widehat\theta(\mathbf y)\mid\mathbb Y\}
		\right)
		\xrightarrow{d}
		\mathcal N\!\left(
		0,\,
		\frac{\mathbf e_0^\top S^{-1}S^*S^{-1}\mathbf e_0}{\sigma^2(\mathbf y) f_{\mathbf Y}(\mathbf y)}
		\right).
		\]
	\end{corollary}
	
	%-------------------------------------------------
	\subsection{MSE and optimal bandwidth}
	
	The conditional mean squared error of $\widehat{D^\alpha\nu}(\mathbf y)$ satisfies
	\(\mathrm{MSE}\{\widehat{D^\alpha\nu}(\mathbf y)\}
	\approx
	h^{2(p+1-|\alpha|)}
	+
	\frac{1}{Nh^{\,s+2|\alpha|}}.
	\)
	Balancing these terms yields the optimal bandwidth rate
	\(h_{\mathrm{opt}}\propto N^{-1/(2p+2+s)},
	\)
	which depends on the polynomial degree $p$ and the covariate dimension $s$,
	but not on the derivative order $|\alpha|$.

	\subsection{Connection to uniform theory}
	
	The pointwise expansions and asymptotic normality results above
	provide the local bias-variance structure underlying the uniform convergence
	rates established in Section~4.
	When combined with entropy bounds for kernel-weighted likelihood classes,
	they justify uniform error rates of the form
	\(h^{\,p+1}
	+
	\sqrt{\frac{\log N}{Nh^{\,s}}},
	\)
	even under data-driven bandwidth selection.
	
		\paragraph{Funding}
	Author state no funding involved.

	\paragraph{Author contributions}
	he author is solely responsible for all aspects of this work, including the conception and design of the study, data generation and analysis, theoretical development, manuscript preparation, and all revisions.

	\paragraph{Conflict of interest}
	The author declare that they have no competing interests.

\end{document}